\documentstyle{amsart}

\newtheorem{Theorem}{Theorem}[section]

\newtheorem{Lemma}[Theorem]{Lemma}

\newtheorem{Corollary}[Theorem]{Corollary}

\newtheorem{Remark}[Theorem]{Remark}

\newtheorem{Definition}[Theorem]{Definition}

\begin{document}

\title{Strong Toroidalization of birational morphisms of 3-folds}

\author{Steven Dale Cutkosky}
\thanks{Research    partially supported by NSF}

\maketitle

\section{Introduction} 

Suppose that $f:X\rightarrow Y$ is a dominant morphism of algebraic varieties, over a field $k$ of characteristic zero.
If $X$ and $Y$ are nonsingular,
$f:X\rightarrow Y$ is toroidal if there are simple normal crossing divisors $D_X$ on $X$ and $D_Y$ on $Y$ such that
$f^{-1}(D_Y)=D_X$, and $f$ is  locally given by monomials in appropriate etale local parameters on $X$.
The precise definition of this concept is in \cite{AK}, \cite{KKMS} and Definition \ref{Def274} of this paper. The problem of toroidalization is to determine, given a dominant morphism $f:X\rightarrow Y$,
 if there exists a commutative diagram 
 \begin{equation}\label{eq393}
 \begin{array}{rcl}
 X_1&\stackrel{f_1}{\rightarrow}&Y_1\\
 \Phi\downarrow&&\downarrow\Psi\\
 X&\stackrel{f}{\rightarrow}&Y
 \end{array}
 \end{equation}
 such that $\Phi$ and $\Psi$ are products of blow ups of nonsingular subvarieties, $X_1$ and $Y_1$ are nonsingular,
 and there exist simple normal crossing divisors $D_{Y_1}$ on $Y_1$ and $D_{X_1}=f^{-1}(D_{Y_1})$ on $X_1$ such that $f_1$ is toroidal
 (with respect to $D_{X_1}$ and $D_{Y_1}$).     This  is stated in Problem 6.2.1. of \cite{AKMW}.
 
 A stronger form of toroidalization is also asked for in Problem 6.2.1 \cite{AKMW}, which we will call strong toroidalization.
 Suppose that $f:X\rightarrow Y$ is a dominant morphism of nonsingular projective varieties , $D_Y$ is a SNC divisor on $Y$ and $D_X=f^{-1}(D_Y)$ is a 
 SNC divisor on $X$ such that the locus $\text{sing}(f)$ where the morphism $f$ is not smooth is contained in $D_X$.
 The problem of strong toroidalization is to determine if there exists a commutative diagram (\ref{eq393})
 such that $\Phi$ and $\Psi$ are products of blow ups of nonsingular centers which are supported in the preimages of
 $D_X$ and $D_Y$ respectively, and make SNCs with the respective preimages of $D_X$ and $D_Y$, and $f_1$ 
 is toroidal with respect to $D_{Y_1}=\Psi^{-1}(D_Y)$ and $D_{X_1}=\Phi^{-1}(D_X)$.

 Toroidalization, and related concepts, have been considered earlier in different contexts,
 mostly for morphisms of surfaces. Strong torodialization is the strongest structure theorem which could be true for
 general morphisms.  The concept of torodialization fails completely in positive characteristic. A simple example is
 shown in \cite{C3}.
 
 In the case when $Y$ is a curve, toroidalization follows from embedded resolution of singularities (\cite{H}).
 When $X$ and $Y$ are surfaces, there are several proofs in print (\cite{AkK}, Corollary 6.2.3 \cite{AKMW}, \cite{CP}, \cite{Mat}).  They all make use of special properties of the
 birational geometry of surfaces.  An outline of proofs of the above cases can be found in the introduction to
 \cite{C3}.

 In \cite{C3}, strong toroidalization  is solved in the case when $X$ is a 3-fold and $Y$ is a surface,
 In Theorem 0.1 of \cite{C5} we prove toroidalization of birational morphisms of 3-folds. 
 In this paper, we prove  strong  toroidalization for birational morphisms of 3-folds.

\begin{Theorem}\label{Theorem3} Suppose that $f:X\rightarrow Y$ is a birational morphism of nonsingular projective 3-folds  over an
algebraically closed field $k$ of characteristic 0. Further suppose that there is a SNC divisor $D_Y$ on $Y$
such that  $D_X=f^{-1}(D_Y)$ is a SNC divisor which contains the singular locus of the map $f$.
 Then there exists a commutative diagram of morphisms
$$
\begin{array}{rll}
X_1&\stackrel{f_1}{\rightarrow}&Y_1\\
\Phi\downarrow&&\downarrow\Psi\\
X&\stackrel{f}{\rightarrow}&Y
\end{array}
$$
where $\Phi,\Psi$ are products of possible blow ups for the preimages of $D_X$ and $D_Y$ respectively,
  and $f_1$ is toroidal with respect to $D_{Y_1}=\Psi^{-1}(D_Y)$ and $D_{X_1}=\Phi^{-1}(D_X)$.
\end{Theorem}

A possible blow up on a nonsingular 3-fold with toroidal structure is the blow up of a point or a nonsingular curve contained in the toroidal structure which makes SNCs with the toroidal structure.

As a consequence of Theorem \ref{Theorem3}, we find the following strong toroidalization theorem for morphisms of (possibly singular) varieties.

\begin{Theorem}\label{Theorem1} Suppose that $f:X\rightarrow Y$ is a birational morphism of 3-folds which are proper over an
algebraically closed field $k$ of characteristic 0. Further suppose that there is an equidimensional codimension 1 reduced subscheme  $D_Y$ of $Y$
such that $D_Y$ contains the singular locus of $Y$, and $D_X=f^{-1}(D_Y)$ contains the singular locus of the map $f$.
 Then there exists a commutative diagram of morphisms
$$
\begin{array}{rll}
X_1&\stackrel{f_1}{\rightarrow}&Y_1\\
\Phi\downarrow&&\downarrow\Psi\\
X&\stackrel{f}{\rightarrow}&Y
\end{array}
$$
where $\Phi,\Psi$ are products of blow ups of nonsingular curves and points supported above $D_X$ and $D_Y$ respectively,
 $D_{Y_1}=\Psi^{-1}(D_Y)$ is a simple normal crossings  divisor on $Y_1$,
$D_{X_1}=f_1^{-1}(D_{Y_1})$ is a simple normal crossings  divisor on $X_1$ and $f_1$ is toroidal with respect to $D_{Y_1}$ and $D_{X_1}$.
\end{Theorem}

The bulk of this paper is devoted to proving the following theorem.

\begin{Theorem}\label{Theorem2}
Suppose that $f:X\rightarrow Y$ is a dominant morphism of nonsingular projective 3-folds
over an algebraically closed field $k$ of characteristic zero, with toroidal structures determined by SNC divisors $D_Y$ on $Y$ and $D_X=f^{-1}(D_Y)$ on $X$
such that  $D_X$ contains the singular locus of $f$.
Then there exists a commutative diagram
$$
\begin{array}{lll}
X_1&\stackrel{f_1}{\rightarrow}&Y_1\\
\Phi\downarrow&&\downarrow\Psi\\
X&\stackrel{f}{\rightarrow}&Y
\end{array}
$$
such that $\Psi$ and $\Phi$ are products  of possible blow ups  for the preimages of $D_Y$, $D_X$ respectively, such that $f_1$ is prepared
for  $D_{Y_1}=\Psi^{-1}(D_Y)$ and $D_{X_1}=\Phi^{-1}(D_X)$, and $D_{X_1}$ is  cuspidal for $f_1$.
\end{Theorem}

The notation used in the statement of Theorem \ref{Theorem2} is defined in Sections 2 and 3 of this paper. From Theorem  \ref{Theorem2} we easily deduce Theorems \ref{Theorem3} and \ref{Theorem1} from results in \cite{C5}.

Theorem \ref{Theorem2} is applicable to arbitrary dominant morphisms of 3-folds, and is a significant step towards a  proof of strong toroidalization of arbitrary dominant morphisms of 3-folds.

If we relax some of the restrictions in the definition of toroidalization, there are other constructions
producing a toroidal morphism $f_1$, which 
are valid for arbitrary dimensions of $X$ and $Y$.
In \cite{AK} it is shown that a diagram (\ref{eq393}) can be constructed where $\Phi$ is weakened to being a
modification (an arbitrary birational morphism).  In \cite{C1}, \cite{C2} and \cite{C4}, it is shown that a diagram (\ref{eq393})
can be constructed where $\Phi$ and $\Psi$ are locally products of blow ups of nonsingular centers and $f_1$ is locally toroidal, but the morphisms $\Phi$, $\Psi$ and $f_1$
may not be separated.  This construction is obtained by patching local solutions, at least one of which contains the center of any given valuation.

It has been shown in \cite{AKMW} and \cite{W2} that weak factorization of birational morphisms holds in characteristic zero,
and arbitrary dimension.  That is, birational morphisms of complete varieties can be factored by an alternating
sequence of blow ups and blow downs of non singular subvarieties.  Weak factorization of birational (toric) morphisms
of toric varieties, (and of birational toroidal morphisms) has been proven by Danilov \cite{D1} and Ewald \cite{E} (for 3-folds), and by Wlodarczyk \cite{W1}, Morelli \cite{Mo} and Abramovich, Matsuki and Rashid \cite{AMR} in general dimensions.

Our Theorem \ref{Theorem3} on strong toroidalization (or the weaker Theorem 0.1 of \cite{C5} on toroidalization), when combined with weak factorization for toroidal morphisms (\cite{AMR}),
gives a new proof of weak factorization of birational morphisms of 3-folds. 
We point out that our proof uses an  analysis of the structure as power series of local germs of a mapping, as opposed to the entirely different proof of weak factorization, using geometric invariant theory, of \cite{AKMW} and \cite{W1}.

The version of weak factorization which we get from Theorem \ref{Theorem3} is  stronger than that obtained in \cite{AKMW}, \cite{W1} or \cite{C5}.

\begin{Corollary}\label{CorWF}
Suppose that $f:X\rightarrow Y$ is a dominant birational morphism of nonsingular projective 3-folds
over an algebraically closed field $k$ of characteristic zero, with toroidal structures determined by SNC divisors $D_Y$ on $Y$ and $D_X=f^{-1}(D_Y)$ on $X$
such that  $D_X$ contains the singular locus of $f$.
Then there exists a commutative diagram of morphisms factoring $f$,
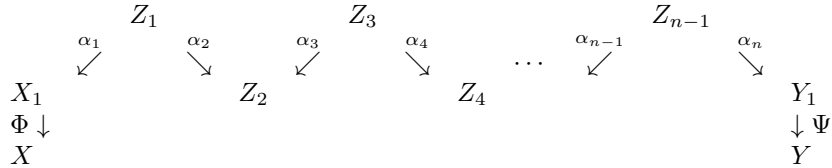
$$
\begin{array}{llllllllllllll}
&&Z_1&&&&Z_3&&&&&Z_{n-1}&&\\
&\stackrel{\alpha_1}{\swarrow}&&\stackrel{\alpha_2}{\searrow}&&\stackrel{\alpha_3}{\swarrow}&&\stackrel{\alpha_4}{\searrow}&&\cdots&
\stackrel{\alpha_{n-1}}{\swarrow}&&\stackrel{\alpha_n}{\searrow}&\\
X_1&&&&Z_2&&&&Z_4&&&&&Y_1\\
\Phi\downarrow&&&&&&&&&&&&&\downarrow\Psi\\
X&&&&&&&&&&&&&Y
\end{array}
$$
such that
\begin{enumerate}
\item[1.] All varieties $X_1$, $Y_1$ and the $Z_i$  are nonsingular, with toroidal structures $D_{X_1}$, $D_{Y_1}$ and $D_{Z_i}$ respectively.
\item[2.] There is a toroidal morphism $f_1:X_1\rightarrow Y_1$ making a strong toroidalization  
$$
\begin{array}{lll}
X_1&\stackrel{f_1}{\rightarrow}& Y_1\\
\downarrow&&\downarrow\\
X&\stackrel{f}{\rightarrow}&Y.\\
\end{array}
$$
\item[3.] The morphisms in the diagram
$$
\begin{array}{llllllllllllll}
&&Z_1&&&&Z_3&&&&&Z_{n-1}&\\
&\stackrel{\alpha_1}{\swarrow}&&\stackrel{\alpha_2}{\searrow}&&\stackrel{\alpha_3}{\swarrow}&&\stackrel{\alpha_4}{\searrow}&&\cdots&
\stackrel{\alpha_{n-1}}{\swarrow}&&\stackrel{\alpha_n}{\searrow}&\\
X_1&&&&Z_2&&&&Z_4&&&&&Y_1\\
\end{array}
$$
are toroidal with respect to their toroidal structures.
\end{enumerate}
\end{Corollary}
 
 The proof of Corollary \ref{CorWF} is immediate from Theorem \ref{Theorem3}, which constructs the commutative diagram 2, and \cite{AMR}, \cite{Mo} or \cite{W1}, which produces the diagram 3.

 The problem of strong factorization, as proposed by Abhyankar \cite{Ab2} and Hironaka \cite{H}, is to factor a birational morphism
$f:X\rightarrow Y$ by constructing a diagram
$$
\begin{array}{lllll}
&&Z&&\\
&\swarrow&&\searrow&\\
X&&\stackrel{f}{\rightarrow}&&Y
\end{array}
$$
where $Z\rightarrow X$ and $Z\rightarrow Y$ factor as products of blow ups of nonsingular subvarieties.
Oda \cite{O} has proposed the analogous problem for (toric) morphisms of toric varieties.

A birational morphism $f:S\rightarrow Y$ of (nonsingular) surfaces can be directly factored by blowing up points
(Zariski \cite{Z1} and Abhyankar \cite{Ab1}), but there are examples showing that a direct factorization is not possible in general
for 3-folds (Shannon \cite {Sh} and Sally\cite{S}).

We also obtain as an immediate corollary of Theorem \ref{Theorem3} the following new result,
which reduces the problem of strong factorization of 3-folds to the case of morphisms of toric varieties
\begin{Corollary}
Suppose that the Oda conjecture on strong factorization of birational  morphisms of 3-dimensional toric varieties is true.
Then the Abhyankar, Hironaka strong factorization conjecture of birational morphisms of complete (characteristic zero) 3-folds is true.
\end{Corollary}

Abhyankar's local factorization conjecture \cite{Ab2}, which is ``strong factorization'' along a valuation,  follows from 
local monomialization (Theorem A \cite{C2}), to reduce to a locally toroidal morphism, and local factorization for
toroidal morphisms along a toroidal valuation Christensen  \cite{Ch} (for 3-folds), and Karu \cite{K} in general dimensions. A proof in the spirit of \cite{Ch} of local factorization of a toroidal morphism in all dimensions, using only elementary properties of determinants, is given in \cite{CS}.

\section{Notation}
Throughout this paper, $k$ will be an algebraically closed field of characteristic zero. A curve, surface or 3-fold is
a quasi-projective variety over $k$ of respective dimension 1, 2 or 3.
If $X$ is a variety, and $p\in X$ is a nonsingular point, then regular parameters at $p$ are regular parameters in ${\cal O}_{X,p}$.
Formal regular parameters at $p$ are regular parameters in $\hat{\cal O}_{X,p}$.
 If $X$ is a variety and $V\subset X$ is a subvariety, then
${\cal I}_V\subset {\cal O}_X$ will denote the ideal sheaf of $V$.
If $V$ and $W$ are subvarieties of a variety $X$, we denote the scheme theoretic intersection 
$Y=\text{spec}({\cal O}_X/{\cal I}_V+{\cal I}_W)$ by $Y=V\cdot W$.

Let $f:X\rightarrow Y$ be a morphism of varieties. We will denote by $\text{sing}(f)$ the closed set of points $p\in X$ such that $f$ is not smooth at $p$. If $D$ is a Cartier divisor on $Y$, then $f^{-1}(D)$ will denote
the reduced divisor $f^*(D)_{red}$.

 Suppose that $a,b,c,d\in{\bf Q}$. Then we will write $(a,b)\le (c,d)$ if
$a\le b$ and $c\le d$.

A toroidal structure on a nonsingular variety $X$ is a simple normal crossing divisor (SNC divisor) $D_X$ on $X$.

We will say that a nonsingular curve $C$ which is a subvariety of a nonsingular 3-fold $X$ with toroidal structure
$D_X$ makes simple normal crossings (SNCs) with $D_X$ if for all $p\in C$, there exist regular parameters
$x,y,z$ at $p$ such that $x=y=0$ are local equations of $C$, and $xyz=0$ contains the support of $D_X$ at $p$.

Suppose that $X$ is a nonsingular 3-fold with toroidal structure $D_X$. If $p\in D_X$ is on the intersection of three components of $D_X$ then $p$ is called a 3-point. If $p\in D_X$
is on the intersection of two components of $D_X$ (and is not a 3-point) then $p$ is called a 2-point. If $p\in D_X$
is not a 2-point or a 3-point, then $p$ is called a 1-point. If $C$ is an irreducible component of the intersection of two
components of $D_X$, then $C$ is called a 2-curve. 

A possible center on a nonsingular 3-fold $X$ with toroidal structure defined by a SNC divisor $D_X$, is a point
on $D_X$ or a nonsingular curve in $D_X$ which makes SNCs with $D_X$. A possible center on a nonsingular surface $S$
with toroidal structure defined by a SNC divisor $D_S$  is a point on $D_S$. We will also call the blow up of a possible center a possible blow up.

Observe that if $\Phi:X_1\rightarrow X$ is the blow up of a possible center, then $D_{X_1}=\Phi^{-1}(D_X)$ is a SNC
divisor on $X_1$. Thus $D_{X_1}$ defines a toroidal structure on $X_1$. All blow ups $\Phi:X_1\rightarrow X$
considered in this paper will be of possible centers, and we will impose the toroidal structure on $X_1$ defined by
$D_{X_1}=\Phi^{-1}(D_X)$.

By a general point $q$ of a variety $V$, we will mean a point $q$ which 
satisfies conditions which hold on some nontrivial open subset of $V$.
The exact open condition which we require will generally be clear from context.
By a general section of a coherent sheaf ${\cal F}$ on a projective variety $X$, we mean the section corresponding 
to a general point of the $k$-linear space $\Gamma(X,{\cal F})$.

If $X$ is a variety, ${k}(X)$ will denote the function field of $X$. A 0-dimensional valuation $\nu$ of ${k}(X)$ is
a valuation of ${k}(X)$ such that $\bold k$ is contained in the valuation ring $V_{\nu}$ of $\nu$ and the residue field of
$V_{\nu}$ is $k$.  If $X$ is a projective variety which is birationally equivalent to $X$, then there exists a
unique (closed) point $p_1\in X_1$ such that $V_{\nu}$ dominates ${\cal O}_{X_1,p_1}$. $p_1$ is called the center of $\nu$
on $X_1$. If $p\in X$ is a (closed) point, then there exists a 0-dimensional valuation $\nu$ of ${k}(X)$ such that $V_{\nu}$ dominates ${\cal O}_{X,p}$
(Theorem 37, Section 16, Chapter VI \cite{ZS}). For $a_1,\ldots, a_n\in {k}(X)$, $\nu(a_1),\ldots, \nu(a_n)$ are rationally dependent if there exist
$\alpha_1,\ldots, \alpha_n\in{\bf Z}$ which are not all zero, such that $\alpha_1\nu(a_1)+\cdots\alpha_n\nu(a_n)=0$ (in the value group of $\nu$).
Otherwise, $\nu(a_1),\ldots, \nu(a_n)$ are rationally independent.

\section{toroidal morphisms and prepared morphisms}

Suppose that $X$ is a nonsingular variety with toroidal structure $D_X$. We will say that an ideal sheaf ${\cal I}\subset {\cal O}_X$ is toroidal
if ${\cal I}$ is locally generated by monomials in local equations of components of $D_X$.

Suppose that $q\in X$. We say that $u,v,w$ are (formal) permissible
parameters at $q$ (for $D_X$) if $u,v,w$ are regular parameters in  $\hat {\cal O}_{X,q}$ such that
\begin{enumerate}
\item[1.] If $q$ is a 1-point, then $u\in{\cal O}_{X,q}$ and $u=0$ is a local equation of
$D_X$ at $q$.
\item[2.] If $q$ is a 2-point then $u,v\in {\cal O}_{X,q}$ and $uv=0$ is a local equation of $D_X$ at $q$.
\item[3.] If $q$ is a 3-point then $u,v,w\in{\cal O}_{X,q}$ and $uvw=0$ is a local equation of
$D_X$ at $q$.
\end{enumerate}
$u,v,w$ are algebraic permissible parameters if we further have that $u,v,w\in{\cal O}_{X,q}$.

\begin{Definition}\label{torf} Let $f:X\rightarrow Y$ be a dominant morphism of nonsingular 3-folds with toroidal structures $D_Y$ on $Y$ and $D_X=f^{-1}(D_Y)$ on $X$ such that $\text{sing}(f)\subset D_X$. Suppose that $u,v,w$ are (possibly formal) permissible parameters
at $q\in Y$. Then
$u,v$ are {\bf toroidal forms} at $p\in f^{-1}(q)$ if there exist permissible parameters $x,y,z$
in $\hat{\cal O}_{X,p}$ such that
\begin{enumerate}
\item[1.]  $q$ is a 2-point or a 3-point, $p$ is a 1-point and 
\begin{equation}\label{eqTF1}
u=x^a,
v=x^b(\alpha+y)
\end{equation}
where $0\ne \alpha\in  k$.
\item[2.] $q$ is 2-point or a 3-point, $p$ is a 2-point and 
\begin{equation}\label{eqTF21}
u=x^ay^b,
v=x^cy^d
\end{equation}
with $ad-bc\ne 0$.
\item[3.] $q$ is a 2-point or a 3-point, $p$ is a 2-point and 
\begin{equation}\label{eqTF22}
u=(x^ay^b)^k,
v=(x^ay^b)^t(\alpha+z)
\end{equation}
where $0\ne\alpha\in  k$, $a,b,k,t>0$ and $\text{gcd}(a,b)=1$.
\item[4.] $q$ is a 2-point or a 3-point, $p$ is a 3-point and 
\begin{equation}\label{eqTF3}
u=x^ay^bz^c,
v=x^dy^ez^f
\end{equation}
where 
$$
\text{rank}\left(\begin{matrix} a&b&c\\ d&e&f\end{matrix}\right)=2.
$$
\item[5.] $q$ is a 1-point, $p$ is a 1-point and 
\begin{equation}\label{eqTF01}
u=x^a,
v=y
\end{equation}
\item[6.] $q$ is a 1-point, $p$ is a 2-point and 
\begin{equation}\label{eqTF02}
u=(x^ay^b)^k,
v=z
\end{equation}
with $a,b,k>0$ and $\text{gcd}(a,b)=1$
\end{enumerate}
\end{Definition}

Regular parameters $x,y,z$ as in Definition \ref{torf} will be called permissible
parameters for $u,v,w$ at $p$.

\begin{Definition}\label{Def274} (\cite{KKMS}, \cite{AK})
A normal variety $\overline X$ with a SNC divisor $D_{\overline X}$ on $\overline X$ is called toroidal if for every point $p\in \overline X$ there exists an affine
toric variety $X_{\sigma}$, a point $p'\in X_{\sigma}$ and an isomorphism of $k$-algebras
$$
\hat{\cal O}_{\overline X,p}\cong \hat{\cal O}_{X_{\sigma},p'}
$$
such that the ideal of $D_{\overline X}$ corresponds to the ideal of $X_{\sigma}-T$ (where $T$ is the  torus in $X_{\sigma}$). Such a pair
$(X_{\sigma},p')$ is called a local model at $p\in \overline X$. $D_{\overline X}$ is called a toroidal structure on $\overline X$.

A dominant morphism $\Phi:\overline X\rightarrow \overline Y$ of toroidal varieties with SNC divisors $D_{\overline Y}$
on $\overline Y$ and  $D_{\overline X}=\Phi^{-1}(D_{\overline Y})$ on $\overline X$, 
is called toroidal at $p\in\overline X$, and we will say that $p$ is a toroidal point of $\Phi$ if with $q=\Phi(p)$, there exist local models
$(X_{\sigma},p')$ at $p$, $(Y_{\tau},q')$ at $q$ and a toric morphism $\Psi:X_{\sigma}\rightarrow Y_{\tau}$ such that the following
diagram commutes:
$$
\begin{array}{rll}
\hat{\cal O}_{\overline X,p}&\leftarrow &\hat{\cal O}_{X_{\sigma},p'}\\
\hat\Phi^*\uparrow&&\hat\Psi^*\uparrow\\
\hat{\cal O}_{\overline Y,q}&\leftarrow&\hat{\cal O}_{Y_{\tau},q'}.
\end{array}
$$
$\Phi:\overline X\rightarrow \overline Y$ is called toroidal (with respect to $D_{\overline Y}$ and $D_{\overline X}$) if $\Phi$ is toroidal at all $p\in \overline X$.
\end{Definition}
The following is the list of toroidal forms for a dominant morphism $f:X\rightarrow Y$
of nonsingular 3-folds with toroidal structure $D_Y$ and $D_X=f^{-1}(D_X)$. Suppose that $p\in D_X$, $q=f(p)\in D_Y$,
and $f$ is toroidal at $p$.
Then there exist permissible parameters $u,v,w$ at $q$ and permissible parameters $x,y,z$ for $u,v,w$ at $p$ such that
one of the following forms hold:
\begin{enumerate}
\item[1.] $p$ is a 3-point and $q$ is a  3-point,
$$
\begin{array}{ll}
u&=x^ay^bz^c\\
v&=x^dy^ez^f\\
w&=x^gy^hz^i,
\end{array}
$$
where $a,b,d,e,f,g,h,i\in{\bf N}$ and
$$
\text{Det}\left(\begin{array}{lll}
a&b&c\\
d&e&f\\
g&h&i
\end{array}\right)\ne 0.
$$
\item[2.] $p$ is a 2-point and $q$ is a 3-point,
$$
\begin{array}{ll}
u&=x^ay^b\\
v&=x^dy^e\\
w&=x^gy^h(z+\alpha)
\end{array}
$$
with $0\ne \alpha\in  k$ and $a,b,d,e,g,h\in{\bf N}$ satisfy $ae-bd\ne 0$.
\item[3.] $p$ is a 1-point and $q$ is a 3-point,
$$
\begin{array}{ll}
u&=x^a\\
v&=x^d(y+\alpha)\\
w&=x^g(z+\beta)
\end{array}
$$
with $0\ne \alpha,\beta\in  k$, $a,d,g>0$.
\item[4.] $p$ is a 2-point and $q$ is a 2-point,
$$
\begin{array}{ll}
u&=x^ay^b\\
v&=x^dy^e\\
w&=z
\end{array}
$$
with $ae-bd\ne 0$
\item[5.] $p$ is a 1-point and $q$ is a 2-point,
$$
\begin{array}{ll}
u&=x^a\\
v&=x^d(y+\alpha)\\
w&=z
\end{array}
$$
with $0\ne \alpha\in  k$, $a,d>0$.
\item[6.] $p$ is a 1-point and $q$ is a 1-point,
$$
\begin{array}{ll}
u&=x^a\\
v&=y\\
w&=z
\end{array}
$$
with $a>0$.
\end{enumerate}

A prepared morphism $\Phi_X:X\rightarrow S$ from a nonsingular 3-fold $X$ to a nonsingular  surface $S$ (with respect to toroidal structures $D_S$ and $D_X=\Phi_X^{-1}(D_S)$) is defined in Definition 6.5 \cite{C3}. A prepared morphism from a 3-fold to a 3-fold is defined
in Definition 2.4 \cite{C5}. This definition assumes that $f:X\rightarrow Y$ is birational, but this definition is perfectly valid for a generically
finite morphism of 3-folds.

\begin{Remark}\label{Remark1}
Suppose that $f:X\rightarrow Y$ is a dominant proper morphism of nonsingular 3-folds with toroidal structures determined by SNC divisors
$D_Y$ on $Y$ and $D_X=f^{-1}(D_Y)$ on $X$, and $D_X$ contains the singular locus of the morphism $f$.
With our assumptions on $f$, $f$ is generically finite.
Recall that the fundamental locus of a generically finite morphism $f:X\rightarrow Y$ of nonsingular proper varieties is
$\{p\in Y\mid \text{dim }f^{-1}(p)>0\}$. The fundamental locus is a closed set of codimension $\ge 2$ in $Y$.
Let $\overline X$ be the normalization of $Y$ in the function field of $X$, with induced finite morphism $\lambda:\overline X\rightarrow Y$. The branch
locus of $\lambda$ is contained in the SNC divisor $D_Y$.  Let $E$ be  an irreducible  component of $D_Y$. By Abhyankar's lemma, the irreducible
components of $\lambda^{-1}(E)$ are disjoint. Thus the irreducible components of $D_X$ which dominate $E$ are disjoint.
\end{Remark}

\begin{Definition}\label{Def31}
A dominant morphism $f:X\rightarrow Y$ of nonsingular 3-folds with toroidal structures determined by SNC divisors $D_Y$ on $Y$, and $D_X=f^{-1}(D_Y)$ on $X$ such that the singular locus of $f$ is contained in $D_X$ is {\bf prepared} for $D_Y$ and $D_X$ if:
\begin{enumerate}
\item[1.] If $q\in Y$ is a 3-point,  $u,v,w$ are  permissible parameters at $q$
and $p\in f^{-1}(q)$, then $u,v$ and $w$ are each a unit (in $\hat{\cal O}_{X,p}$) times a monomial in local equations of the toroidal
structure $D_X$ at $p$ . Furthermore, there exists a permutation of $u,v,w$ such that
$u,v$ are toroidal forms  at $p$.
\item[2.] If $q\in Y$ is a 2-point, $u,v,w$ are permissible parameters at $q$ and $p\in f^{-1}(q)$, then either
\begin{enumerate}
\item $u,v$ are toroidal forms  at $p$ or
\item $p$ is a 1-point and there exist regular parameters $x,y,z\in\hat{\cal O}_{X,p}$ such that there is
 an expression 
$$
\begin{array}{ll}
u&=x^a\\
v&=x^c(\gamma(x,y)+x^dz)\\
w&=y
\end{array}
$$
where $\gamma$ is a unit series and $x=0$ is a local equation of $D_X$, or
\item $p$ is a 2-point and there exist regular parameters $x,y,z$ in $\hat{\cal O}_{X,p}$ such that there is
 an expression 
$$
\begin{array}{ll}
u&=(x^ay^b)^k\\
v&=(x^ay^b)^l(\gamma(x^ay^b,z)+x^cy^d)\\
w&=z
\end{array}
$$
where $a,b>0$, $\text{gcd}(a,b)=1$, $ad-bc\ne 0$, $\gamma$ is a unit series and $xy=0$ is a local equation of $D_X$.
\end{enumerate}
\item[3.] If $q\in Y$ is a 1-point,  and $p\in f^{-1}(q)$,
 then there exist  permissible parameters $u,v,w$ at $q$ such that 
  $u,v$ are
 toroidal forms  at $p$.
\end{enumerate}
\end{Definition}

\begin{Definition}\label{Def247} Suppose that $f:X\rightarrow Y$ is a
prepared morphism of nonsingular 3-folds with toroidal structures $D_Y$ and $D_X=f^{-1}(D_Y)$. Then $D_X$ is cuspidal for $f$ if:
\begin{enumerate}
\item[1.] If $E$ is a component of $D_X$ which does not contain a 3-point then $f$ is toroidal
in a Zariski open neighborhood of $E$.
\item[2.] If $C$ is a 2-curve of $X$ which does not contain a 3-point then $f$ is toroidal
in a Zariski open neighborhood of $C$.
\end{enumerate}
\end{Definition}

\begin{Definition} Suppose that $X$ is a nonsingular 3-fold with toroidal structure determined by a SNC divisor $D_X$.
We will say that $D_X$ is strongly cuspidal if every component of $D_X$ contains a 3-point and every 2-curve of $D_X$ contains a 3-point.
\end{Definition}

\section{Preparation above 2 and 3-points}

\begin{Lemma}\label{Lemma1} Suppose that $f:X\rightarrow Y$ is a dominant morphism of nonsingular projective 3-folds
with toroidal structures determined by SNC divisors $D_Y$ and $D_X=f^{-1}(D_Y)$ such that  $D_X$ contains the singular locus of $f$.
Then there exist a commutative diagram
$$
\begin{array}{lll}
X_1&\stackrel{f_1}{\rightarrow}&Y_1\\
\Phi\downarrow&&\downarrow\Psi\\
X&\stackrel{f}{\rightarrow}&Y
\end{array}
$$
such that $\Phi$ and $\Psi$ are products  of blow ups of 2-curves and $f_1$ is toroidal above all 3-points of $Y_1$.
\end{Lemma}
\begin{pf} Suppose that $\nu$ is a 0-dimensional valuation of $k(X)$. We will say that $\nu$ is resolved for $f$ if the center of $\nu$ on $Y$ is not a 3-point or if
the center of $\nu$ on $Y$ is a 3-point, and $f$ is toroidal at the center of $\nu$ on $X$.

Being resolved is an open condition on the Zariski-Riemann manifold of $X$, and if $\nu$ is resolved for $f$ and
$$
\begin{array}{lll}
X_1&\stackrel{f_1}{\rightarrow}&Y_1\\
\Phi_1\downarrow&&\downarrow\Psi_1\\
X_1&\stackrel{f}{\rightarrow}&Y
\end{array}
$$
is a commutative diagram such that $\Phi_1$ and $\Psi_1$ are products of blow ups of 2-curves, then $\nu$ is resolved for $f_1$.

Suppose that  $\nu$ is a 0-dimensional valuation of $k(X)$ such that  the center $q$ of $\nu$ on $Y$ is a 3-point.
Let $p$ be the center of $\nu$ on $X$. Let $u,v,w$ be permissible parameters at $q$. 
\vskip .2truein
\noindent{\bf Case 1} Suppose that $\nu(u),\nu(v),\nu(w)$ are rationally independent. Since $uvw=0$ is a local equation of $D_X$ at $p$,
there exist   regular  parameters $x,y,z$ in ${\cal O}_{X,p}$ such that
$xyz=0$ contains the germ of $D_X$ at $p$, and
 we have an expression
$$
\begin{array}{ll}
u&=x^ay^bz^c\gamma_1\\
v&=x^dy^ez^f\gamma_2\\
w&=x^gy^hz^i\gamma_3
\end{array}
$$
where the $\gamma_i$ are units in ${\cal O}_{X,p}$.  Since $\nu(u),\nu(v),\nu(w)$ are rationally independent,
$\nu(x),\nu(y),\nu(z)$ are also rationally independent and
$$
\text{Det}\left(\begin{array}{lll}
a&b&c\\
d&e&f\\
g&h&i
\end{array}\right)\ne 0
$$
which implies that $p$ is a 3-point and $f$ is toroidal at $p$. Thus $\nu$ is resolved for $f$.
\vskip .2truein
\noindent{\bf Case 2} Suppose that $\nu(u),\nu(v)$ are rationally dependent. After possibly interchanging $u,v,w$ we reduce to this case.
Let $C$ be the 2-curve of $Y$ with local equation $u=v=0$ at $q$.
There exists a sequence of blow ups of 2-curves $\Psi_{\nu}:Y_{\nu}\rightarrow Y$ such that the center of $\nu$ on $Y_{\nu}$ is not a 3-point.

$Y_{\nu}$ is the blow up of a toroidal ideal sheaf ${\cal J}_{\nu}$ of ${\cal O}_Y$. Since $f^{-1}(D_Y) =D_X$, ${\cal J}_{\nu}{\cal O}_X$
is also a toroidal ideal sheaf.
By Lemma 2.11 \cite{C5}, there exists a sequence of blow ups of 2-curves $\Phi_{\nu}:X_{\nu}\rightarrow X$ such that there is a commutative diagram of morphisms
$$
\begin{array}{lll}
X_{\nu}&\stackrel{f_{\nu}}{\rightarrow}&Y_{\nu}\\
\Phi_{\nu}\downarrow&&\downarrow\Psi_{\nu}\\
X&\stackrel{f}{\rightarrow}&Y.
\end{array}
$$

Thus $\nu$ is resolved for $f_{\nu}$.

It follows from compactness of the Zariski Riemann manifold of $X$ \cite{Z}, that there exists a positive integer $n$ and commutative diagrams
$$
\begin{array}{lll}
X_i&\stackrel{f_i}{\rightarrow}&Y_i\\
\Phi_i\downarrow&&\downarrow\Psi_i\\
X&\stackrel{f}{\rightarrow}&Y
\end{array}
$$
for $1\le i\le n$ such that $\Phi_i$ and $\Psi_i$ are products of blow ups of 2-curves, and every 0-dimensional valuation $\nu$ of $k(X)$ is resolved for some
$f_i$.

$Y_i$ is the blow up of a toroidal ideal sheaf ${\cal J}_i$ of ${\cal O}_Y$ and $X_i$ is the blow up of a toroidal ideal sheaf ${\cal I}_i$ of ${\cal O}_X$. 
Thus there exists a sequence of blow ups of 2-curves $Y^*\rightarrow Y$ such that $\prod_i{\cal J}_i{\cal O}_{Y^*}$ is invertible,
by Lemma 2.11 \cite{C5}.
$Y^*$ is thus the blow up of a toroidal ideal sheaf ${\cal J}\subset {\cal O}_Y$, so that ${\cal J}{\cal O}_X$ is also a toroidal
ideal sheaf. By Lemma 2.11 \cite{C5}, there exists a sequence of blow ups of 2-curves
$X^*\rightarrow X$ such that ${\cal J}\prod{\cal I}_i{\cal O}_{X^*}$ is invertible.
Thus for $1\le i\le n$ there exist commutative diagrams of morphisms
$$
\begin{array}{lll}
X^*&\stackrel{f^*}{\rightarrow}&Y^*\\
\Phi_i^*\downarrow&&\downarrow\Psi_i^*\\
X_i&\stackrel{f_i}{\rightarrow}&Y_i\\
\downarrow&&\downarrow\\
X&\rightarrow&Y.
\end{array}
$$

Suppose that $\nu$ is a 0-dimensional valuation of $k(X)$. If the center of $\nu$ on $Y^*$ is a 3-point, then the center of $\nu$ on $Y_i$
is a 3-point for all $i$, since $\Psi_i^*$ is toroidal. There exists an $i$ such that $\nu$ is resolved for $f_i$. Thus $f_i$ is toroidal at the center of $\nu$ on $X_i$.
Since $\Phi_i^*$ and $\Psi_i^*$ are toroidal, $f^*$ is toroidal at the center of $\nu$. Thus $\nu$ is resolved for $f^*$.
Since all 0-dimensional valuations of $k(X)$ are resolved for $f^*$, it follows that $f^*$ is toroidal above all 3-points of $Y^*$, and we have achieved the
conclusions of the lemma.
\end{pf}

\begin{Lemma}\label{Lemma2} Suppose that $f:X\rightarrow Y$ is a dominant morphism of nonsingular projective 3-folds, 
with toroidal structures determined by  SNC divisors $D_Y$ and $D_X=f^{-1}(D_Y)$ such that $D_X$ contains the singular locus of $f$.
Further suppose that $f$ is toroidal above all 3-points of $Y$.
Then there exists a commutative diagram
$$
\begin{array}{lll}
X_1&\stackrel{f_1}{\rightarrow}&Y_1\\
\Phi\downarrow&&\downarrow\Psi\\
X&\stackrel{f}{\rightarrow}&Y
\end{array}
$$
such that $\Psi$ and $\Psi$ are products  of blow ups of 2-curves, $f_1$ is toroidal above all 3-points of $Y_1$,
and $f_1$ is prepared (and satisfies 2a) of Definition \ref{Def31}) above all 2-points of $Y_1$.
\end{Lemma}
\begin{pf} Suppose that $\nu$ is a 0-dimensional valuation of $k(X)$. We will say that $\nu$ is resolved for $f$ if the center of $\nu$ on $Y$ is
a 1-point or if the center of $\nu$ on $Y$ is a 2-point and $f$ is prepared at the center of $\nu$ on $X$
(and satisfies 2a) of Definition \ref{Def31}), or if the center of $\nu$ on $Y$
is a 3-point, and $f$ is toroidal at the center of $\nu$ on $X$. 

Being  resolved is an open condition on the Zariski-Riemann manifold of $X$.
Suppose that 
$$
\begin{array}{lll}
X_1&\stackrel{f_1}{\rightarrow}&Y_1\\
\Phi\downarrow&&\downarrow\Psi\\
X&\stackrel{f}{\rightarrow}&Y
\end{array}
$$
is a commutative diagram of morphisms such that $\Phi$ and $\Psi$ are products of blow ups of 2-curves. If $\nu$ is a 0-dimensional valuation of $k(X)$ such that $\nu$ is resolved for $f$, then $\nu$ is resolved for $f_1$.

Suppose that $q\in Y$ is a 2-point, and $\nu$ is a 0-dimensional valuation of $k(X)$ such that $q$ is the center of $\nu$ on $Y$.
Let $p$ be the center of $\nu$ on $X$. Let $u,v,w$ be permissible parameters at $q$, so that $u=v=0$ are local equations of the 2-curve $C$ through $q$. 
\vskip .2truein
\noindent{\bf Case 1} Suppose that $\nu(u),\nu(v)$ are rationally independent. Since $uv=0$ is a local equation of $D_X$ at $p$,
there exist   regular  parameters $x,y,z$ in ${\cal O}_{X,p}$ such that
$xyz=0$ contains the germ of $D_X$ in ${\cal O}_{X,p}$, and
 we have an expression
$$
\begin{array}{ll}
u&=x^ay^bz^c\gamma_1\\
v&=x^dy^ez^f\gamma_2\\
\end{array}
$$
where the $\gamma_i$ are units in ${\cal O}_{X,p}$. Since $\nu(u),\nu(v)$ are rationally independent, 
$$
\text{rank}\left(\begin{array}{lll}
a&b&c\\
d&e&f\\
\end{array}\right)=2
$$
which implies that $p$ is a 2 or 3-point and $f$ is prepared at $p$ (and satisfies 2a) of Definition \ref{Def31}). Thus $\nu$ is resolved for $f$.
\vskip .2truein
\noindent{\bf Case 2} Suppose that $\nu(u),\nu(v)$ are rationally dependent. 
There exists a sequence of blow ups of 2-curves $\Psi_{\nu}:Y_{\nu}\rightarrow Y$ such that the center of $\nu$ on $Y_{\nu}$ is a 1-point.

$Y_{\nu}$ is the blow up of a toroidal ideal sheaf ${\cal J}_{\nu}$ of ${\cal O}_Y$. Since $f^{-1}(D_Y)=D_X$, ${\cal J}_{\nu}$ is also
a toroidal ideal sheaf.
By Lemma 2.11 \cite{C5}, and induction on the number of 2-curves blown up by $\Psi_{\nu}$, there exists a sequence of blow ups of 2-curves $\Phi_{\nu}:X_{\nu}\rightarrow X$ such that there is a commutative diagram of morphisms
$$
\begin{array}{lll}
X_{\nu}&\stackrel{f_{\nu}}{\rightarrow}&Y_{\nu}\\
\Phi_{\nu}\downarrow&&\downarrow\Psi_{\nu}\\
X&\stackrel{f}{\rightarrow}&Y.
\end{array}
$$

Thus $\nu$ is resolved for $f_{\nu}$.

It follows from compactness of the Zariski Riemann manifold of $X$ \cite{Z} that there exists a positive integer $n$ and commutative diagrams
$$
\begin{array}{lll}
X_i&\stackrel{f_i}{\rightarrow}&Y_i\\
\Phi_i\downarrow&&\downarrow\Psi_i\\
X&\stackrel{f}{\rightarrow}&Y
\end{array}
$$
for $1\le i\le n$ such that $\Phi_i$ and $\Psi_i$ are products of blow ups of 2-curves, and every valuation $\nu$ of $k(X)$ is resolved for some
$f_i$.

$Y_i$ is the blow up of a toroidal ideal sheaf ${\cal J}_i$ of ${\cal O}_Y$ and $X_i$ is the blow up of a toroidal ideal sheaf ${\cal I}_i$
of ${\cal O}_X$. 
Thus there exists a sequence of blow ups of 2-curves $Y^*\rightarrow Y$ such that $\prod_i{\cal J}_i{\cal O}_{Y^*}$ is invertible,
by Lemma 2.11 \cite{C5}.
$Y^*$ is thus the blow up of a toroidal ideal sheaf ${\cal J}\subset {\cal O}_Y$. Thus ${\cal J}{\cal O}_X$ is also a toroidal ideal sheaf. By Lemma 2.11 \cite{C5}, there exists a sequence of blow ups of 2-curves
$X^*\rightarrow X$ such that ${\cal J}\prod{\cal I}_i{\cal O}_{X^*}$ is invertible.
Thus for $1\le i\le n$, there exist commutative diagrams of morphisms
$$
\begin{array}{lll}
X^*&\stackrel{f^*}{\rightarrow}&Y^*\\
\Phi_i^*\downarrow&&\downarrow\Psi_i^*\\
X_i&\stackrel{f_i}{\rightarrow}&Y_i\\
\downarrow&&\downarrow\\
X&\rightarrow&Y.
\end{array}
$$
Since $\Phi_i^*$ and $\Psi_i^*$ are the blow ups of toroidal ideal sheaves, they are toroidal morphisms.

Suppose that $\nu$ is a 0-dimensional valuation of $k(X)$. If the center of $\nu$ on $Y$ is a 3-point then $f^*$ is resolved at the center of $\nu$ on $X^*$. 
In particular, if the center of $\nu$ on $Y^*$ is a 3-point, then the center of $\nu$ on $Y$ is a 3-point and $\nu$ is resolved for $f^*$.
Suppose that the center of $\nu$ on $Y^*$ is 2-point, and the center of $\nu$ on $Y$ is not a 3-point. Then the center of
$\nu$ on $Y_i$
is a 2-point for all $i$. There exists an $i$ such that $\nu$ is resolved for $f_i$. Thus $f_i$ is prepared (and satisfies 2a) of Definition \ref{Def31}) at the center of $\nu$ on $X_i$.
Since $\Phi_i^*$ and $\Psi_i^*$ are toroidal, $f^*$ is prepared (and satisfies 2a) of Definition \ref{Def31}) at the center of $\nu$. Thus $\nu$ is resolved for $f^*$.
Since all 0-dimensional valuations of $k(X)$ are resolved for $f^*$, it follows that $f^*$ is toroidal above all 3-points of $Y^*$, and prepared above all 2-point of
$Y^*$, and we have achieved the
conclusions of the lemma.
\end{pf}

\begin{Lemma}\label{Lemma3}  Suppose that $f:X\rightarrow Y$ satisfies the conclusions of Lemma \ref{Lemma2}. Suppose that  $H$ is a general
hyperplane section of $Y$. Then $f$ is prepared above all points of $H$.
\end{Lemma}

\begin{pf} Bertini's theorem implies that $H$ is nonsingular and makes SNCs with $D_Y$. Further, $H'=f^{-1}(H)$ is nonsingular and makes SNCs with
$D_X$. Thus $H$ contains no 3-points of $Y$ and $H'$ contains no 3-points.

Suppose that $q\in H\cap D_Y$ is a 1-point, and $p\in f^{-1}(q)$.  Let $u,v,w$ be regular parameters in ${\cal O}_{Y,q}$
such that $u=0$ is a local equation of $D_Y$ at $q$, and $w=0$ is a local equation of $H$. Then we have regular paramaters $x,y,z$ in ${\cal O}_{X,p}$
such that either $p$ is a 1 point with $x=0$ a local equation of $D_X$ or $p$ is a 2-point with $xy=0$  a local equation of $D_X$ at $p$.
We then have an expression $u=x^a\gamma, w=z$ or
$u=x^ay^b\gamma, w=z$  where $\gamma$ is a unit in ${\cal O}_{X,p}$. Thus $f$ is prepared at $p$.
\end{pf}

\begin{Corollary}\label{Cor3} Suppose that $f:X\rightarrow Y$ satisfies the conclusions of Lemma \ref{Lemma2}. Then there exists a finite set of 1-points
$\Omega\subset Y$ such that $f$ is prepared above $Y-\Omega$.
\end{Corollary}
\begin{pf} The locus of points in $X$ where $f$ is prepared is an open set. Since $f$ is proper, the image $\Omega$ of the closed set of points where $f$ 
is not prepared is closed in $Y$. Since a general hyperplane section of $Y$ is disjoint from $\Omega$ by Lemma \ref{Lemma3}, $\Omega$ must be a finite
set of points.
\end{pf}

\begin{Lemma}\label{Lemma4} Suppose that $f:X\rightarrow Y$ is a proper dominant morphism of nonsingular 3-folds and 
$\pi:Y\rightarrow S$ is a smooth dominant morphism onto a nonsingular surface $S$. Let $g=\pi\circ f$. Suppose that $C$ is a nonsingular curve of $S$,
$D=\pi^{-1}(C)$ and $D'=f^{-1}(D)$. Suppose that $D'$ is a SNC divisor on $X$ which contains the singular locus of $g$ and the singular locus of $f$. Suppose that
$\overline q\in C$ is a point, and that $g$ is toroidal and prepared (with respect to $C$ and $D'$) away from points above finitely many points $\Omega=\{q_1,\ldots, q_m\}\subset\gamma=\pi^{-1}(\overline q)$.  Further suppose that $f$ is finite above a general point of $\gamma$.
Then there exists a commutative diagram of morphisms
$$
\begin{array}{lll}
X_1&\stackrel{\Phi_1}{\rightarrow}&X\\
g_1&\searrow&\downarrow g\\
&&S
\end{array}
$$
such that $\Phi_1$ is a product of possible blow ups  for the preimage of $D'$ supported above $\Omega$ and $g_1$ is prepared (with respect to $C$ and $\Phi_1^{-1}(D')$)
in a neighborhood
of all components $F$ of $\Phi_1^{-1}(D')$ which dominate  $D$  and in a neighborhood of all components $F$ of $\Phi_1^{-1}(D')$ which dominate a curve of $Y$.
\end{Lemma}

\begin{pf} Let $u,w$ be regular parameters in ${\cal O}_{S,\overline q}$ such that $u=0$ is a local equation of $C$ at $\overline q$.
Let $C'$ be the curve on $S$ with local equation $w=0$ at $\overline q$. Let $A=\pi^{-1}(C')$.

Since it suffices to prove the lemma above a neighborhood of $\overline q$ in $S$, we may assume that $E=C+C'$ is a SNC divisor  on $S$ whose only
singular point is $\overline q$.
 Since $g$ is toroidal  away from points above $\Omega$, we have that $g^{-1}(E)$ defines a SNC divisor on $X$
away from points above $\Omega$.
There exists a morphism $\Phi_1:X_1\rightarrow X$ which is a sequence of possible blow ups for the preimage of $D'$ supported above $\Omega$ such that
with $g_1=g\circ\Phi_1:X_1\rightarrow S$, $g_1^{-1}(E)$ is a SNC divisor, and $(f\circ\Phi_1)^{-1}(q_i)$ are divisors for all $q_i\in\Omega$.
We may further assume that the union $\overline A$ of codimension 1 subvarieties of $X_1$ which dominate $A$ are disjoint, since they are disjoint
away from the preimage of $\Omega$.

Let $\overline D$ be the union of codimension 1 subvarieties of $X_1$ which dominate $D$, so that $\overline D$ is a 
disjoint union of irreducible components of $D''=g_1^{-1}(C)$ (by Remark \ref{Remark1}). 

Suppose that $p\in\overline D$ and $f\circ\Phi_1(p)=q_i\in\Omega$. Then $p$ must be a 2-point or a 3-point.
We have regular parameters $x,y,z$ in $\hat{\cal O}_{X_1,p}$ such that one of the following cases hold:
\begin{enumerate}
\item[1.] $p$ is a 2-point and
$$
u=x^ay^b, w=y^c
$$
where $x=0$ is a local equation of $\overline D$, $u=0$ is a local equation of $D''$ and $a,b>0$.
\item[2.] $p$ is a 2-point, 
$$
u=x^ay^b, w=y^cz
$$
where $x=0$ is a local equation of $\overline D$, $u=0$ is a local equation of $D''$,  $a,b,c>0$ and $z=0$ is a local equation of $\overline A$.
\item[3.] $p$ is a 3-point and
$$
u=x^ay^bz^c, w=y^dz^e
$$
where $x=0$ is a local equation of $\overline D$, $u=0$ is a local equation of $D''$ and $a,b,c,d,e>0$.
\end{enumerate}
Thus $g_1$ is prepared in a neighborhood of $\overline D$.

Now suppose that  $F$ is a component of $D''$ which dominates a curve of $Y$ and $p\in F$ is such that
$f\circ\Phi_1(p)=q_i\in\Omega$. Then $p$ must be a 2-point or a 3-point. By our assumption that $f$ is finite above a general point of $\gamma$, 
$F$ dominates the curve $C$ of $S$. Thus we have regular parameters $x,y,z$ in $\hat{\cal O}_{X_1,p}$ such that one of the following cases hold:
\begin{enumerate}
\item[1.] $p$ is a 2-point and
$$
u=x^ay^b, w=y^c
$$
where $x=0$ is a local equation of $F$, $u=0$ is a local equation of $D''$ and $a,b>0$.
\item[2.] $p$ is a 2-point, 
$$
u=x^ay^b, w=y^cz
$$
where $x=0$ is a local equation of $F$, $u=0$ is a local equation of $D''$, $a,b>0$ and $z=0$ is a local equation of  $\overline A$.
\item[3.] $p$ is a 3-point and
$$
u=x^ay^bz^c, w=y^dz^e
$$
where $x=0$ is a local equation of $F$, $u=0$ is a local equation of $D''$ and $a,b,c>0$.
\end{enumerate}
Thus $g_1$ is prepared in a neighborhood of $F$.
\end{pf}

\begin{Lemma}\label{LemmaA} Suppose that $f:X\rightarrow Y$ is a dominant morphism of nonsingular 3-folds
with toroidal structures determined by SNC divisors $D_Y$ and $D_X=f^{-1}(D_Y)$ such that  $D_X$ contains the singular locus of $f$.
Further suppose that $f:X\rightarrow Y$ is toroidal and $q\in Y$ is a 2-point. Let $\Psi:Y_1\rightarrow Y$ be the blow
up of $q$. Then there exists a commutative
diagram of morphisms
$$
\begin{array}{lll}
X_1&\stackrel{f_1}{\rightarrow}&Y_1\\
\Phi\downarrow&&\downarrow\Psi\\
X&\stackrel{f}{\rightarrow}&Y
\end{array}
$$
such that $\Phi$ is a sequence of possible blow ups for the preimage of $D_X$ supported above $q$ and $f_1$ is toroidal with respect to $D_{Y_1}=\Psi^{-1}(D_Y)$
and $D_{X_1}=\Phi^{-1}(D_X)$.
\end{Lemma}

\begin{pf} There exist permissible parameters $u,v,w$ at $q$ such that if $p\in f^{-1}(q)$ then there exist permissible parameters $x,y,z$
for $u,v,w$ such that if $p$ is a 1-point, then we have a form 
\begin{equation}\label{eqA1}
u=x^a,v=x^b(\alpha+y), w=z
\end{equation}
with $0\ne\alpha\in k$, and if $p$ is a 2-point, 
\begin{equation}\label{eqA2}
u=x^ay^b, v=x^cy^d, w=z,
\end{equation}
with $ad-bc\ne 0$.
We first show that there exists a sequence of possible blow ups 
\begin{equation}\label{eqA21}
X_m\stackrel{\Phi_m}{\rightarrow}X_{m-1}\rightarrow \cdots\rightarrow X_1\stackrel{\Phi_1}{\rightarrow} X
\end{equation}
obtained by  blow ups of possible centers supported above $q$ such that the rational map $X_m\rightarrow Y_1$ is toroidal wherever it is defined, and if
${\cal I}_q{\cal O}_{X_m,p}$ is not invertible, then there exist regular parameters $x,y,z$ in $\hat{\cal O}_{X_m,p}$  such that one of the following forms
hold: 

$p$ is a 1-point 
\begin{equation}\label{eqA3}
u=x^a,v=x^b(\alpha+y), w=x^cz
\end{equation}
with $\alpha\ne 0$, and $c=0$ or 1, or $p$ is a 2-point 
\begin{equation}\label{eqA4}
u=x^ay^b, v=x^cy^d, w=xz
\end{equation}
with $a,c\ge 1$, or $p$ is a 2-point 
\begin{equation}\label{eqA5}
u=x^ay^b, v=x^cy^d, w=xyz
\end{equation}
with $a,c\ge 1$ and $b,d\ge 1$.

The points $p\in f^{-1}(p)$ such that $u,v,w$ do  not have a form (\ref{eqA3}), (\ref{eqA4}) or (\ref{eqA5}) at $p$ are 2-points of one of the following
forms: 
\begin{equation}\label{eqA6}
u=x^a, v=y^b,w=z,
\end{equation}
in which case $V(x,y,z)$ is the locus in $\text{spec}(\hat{\cal O}_{X,p})$ where ${\cal I}_q\hat{\cal O}_{X,p}$ is not invertible,
\begin{equation}\label{eqA7}
u=x^a, v=x^by^c,w=z
\end{equation}
with $b,c>0$,
in which case $V(x,z)$ is the locus in $\text{spec}(\hat{\cal O}_{X,p})$ where ${\cal I}_q\hat{\cal O}_{X,p}$ is not invertible,
\begin{equation}\label{eqA8}
u=x^ay^b, v=x^cy^d,w=z
\end{equation}
with $a,b,c,d>0$
in which case $V(x,z)\cup V(y,z)$ is the locus in $\text{spec}(\hat{\cal O}_{X,p})$ where ${\cal I}_q\hat{\cal O}_{X,p}$ is not invertible,

Let $Z$ be the closed locus  of points $r$ in $X$ such that ${\cal I}_q{\cal O}_{X,r}$ is not invertible.
The isolated points $p$ in $Z$   have a form (\ref{eqA6}). If $p$ is a non isolated point in $Z$ which is a 2-point,
 then $p$ has a form (\ref{eqA7}) or (\ref{eqA8}).

Suppose that $E $ is a curve in $Z$ such that $E$ contains a 2-point $p$ satisfying (\ref{eqA7}) or (\ref{eqA8}). Then a generic point of $E$
satisfies (\ref{eqA1}) and all 2-points of $E$ must have a form (\ref{eqA7}) or (\ref{eqA8}).

Let $\Phi_1:X_1\rightarrow X$ be the blow up of the finitely many points $p\in X$ of the form (\ref{eqA6}).
Suppose that $p\in X$ is such a point, 
and $p_1\in\Phi_1^{-1}(p)$. Without loss of generality, we may assume that $a\le b$ in (\ref{eqA6}). There are regular parameters
$x_1,y_1,z_1$ in $\hat{\cal O}_{X_1,p_1}$ of one of the following forms: 
\begin{equation}\label{eqA9}
x=x_1, y=x_1(y_1+\alpha), z=x_1(z_1+\beta)
\end{equation}
with $\alpha,\beta\in k$, 
\begin{equation}\label{eqA10}
x=x_1y_1,y=y_1,z=y_1(z_1+\alpha)
\end{equation}
with $\alpha\in k$ or 
\begin{equation}\label{eqA11}
x=x_1z_1, y=y_1z_1, z=z_1.
\end{equation}
Suppose that (\ref{eqA9}) holds. Then $u,v,w$ have a form
$$
u=x_1^a, v=x_1^b(y_1+\alpha)^b, w=x_1(z_1+\beta)
$$
at $p_1$. If $a=1$, then $f\circ \Phi_1$ factors through $Y_1$ at $p_1$ and we have one of the following toroidal forms:

1-point maps to 2-point:
$$
u_1=u=x_1,
v_1=\frac{v}{u}=x_1^{b-1}(y_1+\alpha)^b,
w_1=\frac{w}{u}-\beta=z_1
$$
if $b>a=1$ and $\alpha\ne 0$,

1-point maps to 1-point:
$$
u_1=u=x_1,
v_1=\frac{v}{u}-\alpha= y_1,
w_1=\frac{w}{u}-\beta=z_1
$$
if $b=a=1$, $\alpha\ne 0$,

2-point maps to 2-point:
$$
u_1=u=x_1,
v_1=\frac{v}{u}=x_1^{b-1}y_1^b,
w_1=\frac{w}{u}-\beta=z_1
$$
if $a=1$ and $\alpha=0$.

Suppose that (\ref{eqA9}) holds and $a>1$. If $\beta\ne 0$, we have that $\Phi_1\circ f$ factors through $Y_1$ at $p_1$ and we have a toroidal form,
obtained from a change of variable in
$$
u_1=\frac{u}{w}=x_1^{a-1}(z_1+\beta)^{-1},
v_1=\frac{v}{w}=x_1^{b-1}(z_1+\beta)^{-1}(y_1+\alpha)^b,
w_1=w=x_1(z_1+\beta)
$$
where $p_1$ is 1-point mapping to a 3-point if $\alpha\ne 0$ and $p_1$ is 2-point mapping to a 3-point if $\alpha=0$.

If $\beta=0$ (and $a>1$) then we have
$$
u=x_1^a,
v=x_1^b(y_1+\alpha),
w=x_1z_1
$$
of the form (\ref{eqA3}) if $\alpha\ne 0$ and of the form (\ref{eqA4}) if $\alpha=0$.

Suppose that (\ref{eqA10}) holds. Then at $p_1$, $u,v,w$ have a form:
$$
u=x_1^ay_1^a,
v=y_1^b,
w=y_1(z_1+\alpha).
$$

Assume $b=1$ (which implies $a=1$). then $f\circ\Phi_1$ factors through $Y_1$ at $p_1$, and there is a toroidal form:
$$
u_1=\frac{u}{v}=x_1,
v_1=v=y_1,
w_1=\frac{w}{v}-\alpha=z_1
$$
where $p_1$ is 2-point mapping to a 2-point.

Assume that $b>1$ and $\alpha\ne 0$. Then $f\circ\Phi_1$ factors through $Y_1$ at $p_1$, and there is a toroidal form, obtained from a change
of variable in
$$
u_1=\frac{u}{w}=x_1^{a}y_1^{a-1}(z_1+\alpha)^{-1},
v_1=\frac{v}{w}=y_1^{b-1}(z_1+\alpha)^{-1},
w_1=w=y_1(z_1+\alpha)
$$
where $p_1$ is a 2-point mapping to a 3-point.

If $b>1$ and $\alpha=0$, then we have a form:
$$
u=x_1^ay_1^a,
v=y_1^b,
w=y_1z_1
$$
of the form (\ref{eqA4}).

Suppose that (\ref{eqA11}) holds. Then $p_1$ is a 3-point and $u,v,w$ have a form
$$
u=x_1^az_1^a,
v=y_1^bz_1^b,
w=z_1.
$$
Thus $\Phi_1\circ f$ factors through $Y_1$ at $p_1$ by
$$
u_1=\frac{u}{w}=x_1^az_1^{a-1},
v_1=\frac{v}{w}=y_1^bz_1^{b-1},
w_1=w=z_1,
$$
where $p_1$ is a 3-point mapping to a 3-point.

We have thus completed the analysis of $\Phi_1$.

We now construct (\ref{eqA21}) by induction. Each $X_i$ will be such that the rational map $X_i\rightarrow Y_1$ is toroidal wherever it is defined, and if
$p\in X_i$ is a 2-point such that ${\cal I}_q{\cal O}_{X_1,p}$ is not invertible, then there exist regular parameters $x,y,z$ at $p$ such that $u,v,w$ have one of the forms (\ref{eqA3}),
(\ref{eqA4}), (\ref{eqA5}), (\ref{eqA7}) or (\ref{eqA8}) at $p$. If a form (\ref{eqA7}) or (\ref{eqA8})  holds at $p$, then $\Phi_1\circ\cdots\circ\Phi_i$ is an isomorphism
near $p$.

Each $\Phi_{i+1}:X_{i+1}\rightarrow X_i$ for $i\ge 1$ will be the blow up of a  curve $E_i$ which is a possible center and is the strict transform of a component
of $Z\subset X$.

Suppose that we have constructed (\ref{eqA21}) out to $X_i$, and $p\in X_i$ is a 2-point such that ${\cal I}_q{\cal O}_{X_i,p}$ is not invertible, and
$u,v,w$ do not have a form (\ref{eqA3}), (\ref{eqA4}) or (\ref{eqA5}) at $p$.
Then $u,v,w$ have a form (\ref{eqA7}) or (\ref{eqA8}) at $p$. Let $E=E_i$ be a curve in the locus where ${\cal I}_q{\cal O}_{X_i}$
is not invertible which contains $p$. Let $F$ be the component of $D_{X_i}$ containing $E_i$ We necessarily have $\text{ord}_Fw=0$ and $\text{ord}_Fu>0$,
$\text{ord}_Fv>0$.
Further, $\Phi_1\circ\cdots\circ\Phi_i$ is an isomorphism near $p$. Thus $E$ is the strict transform of a component of $Z$.

Suppose that $p'\in E_i$ is another 2-point. Then at $p'$, since $\text{ord}_Fw=0$, $u,v,w$   must have a form (\ref{eqA7}), (\ref{eqA8}) or
(\ref{eqA4}), where in this last case, $y=z=0$ is a local equation of $E$ and $b,d\ge 1$ (since $\text{ord}_Fw=0$,  $\text{ord}_Fu>0$ and $\text{ord}_Fv>0$). If $p'\in E_i$ is a 1-point, then
$u,v,w$ have a form (\ref{eqA1}) at $p'$, since $\text{ord}_Fw=0$.

 Let $\Phi_{i+1}:X_{i+1}\rightarrow X_i$ be the blow up of $E$
 and $\overline\Phi_{i+1}=\Phi_1\circ\cdots\circ\Phi_{i+1}$.
 
 Suppose that $p\in E$ is a 1-point and $p_1\in\Phi_{i+1}^{-1}(p)$. Then $f\circ\overline\Phi_{i+1}$ is toroidal whenever it is defined,
 and points above $p$ where $f\circ\overline\Phi_{i+1}$ does not factor through $Y_i$ have a form (\ref{eqA3}).
 A detailed analysis of a case including this one is given later in the proof, after (\ref{eqA17}).

Suppose that $p\in E$ is a 2-point of the form (\ref{eqA8}) and $p_1\in \Phi_{i+1}^{-1}(p)$.

There are regular parameters $x_1,y_1,z_1$ in $\hat{\cal O}_{X_i,p_1}$ of one of the following forms: 
\begin{equation}\label{eqA12}
x=x_1, z=x_1(z_1+\alpha)
\end{equation}
with $\alpha\in k$ or 
\begin{equation}\label{eqA13}
x=x_1z_1, z=z_1.
\end{equation}
Suppose that (\ref{eqA12}) holds. We have that $p_1$ is a 2-point, and
$$
u=x_1^ay^b,
v=x_1^cy^d,
w=x_1(z_1+\alpha).
$$
If $\alpha\ne0$, we have that $f\circ\overline \Phi_{i+1}$ factors through $Y_1$ at $p_1$. We have a form:
$$
u_1=\frac{u}{w}=x_1^{a-1}y^b(z_1+\alpha)^{-1},
v_1=\frac{v}{w}=x_1^{c-1}y^d(z_1+\alpha)^{-1},
w_1=x_1(z_1+\alpha)
$$
at the 2-point $p_1$, which maps to a 3-point, and thus is toroidal, after a change of variables.

If $\alpha=0$ in (\ref{eqA12}), we have
$$
u=x_1^ay^b,
v=x_1^cy^d,
w=x_1z_1
$$
of the form (\ref{eqA4}).

If (\ref{eqA13}) holds, then $p_1$ is a 3-point and
$$
u=x_1^ay^bz_1^a,
v=x_1^cy^dz_1^c,
w=z_1.
$$
Thus $f\circ\overline \Phi_{i+1}$ factors through $Y_1$ at $p_1$, and we have a toroidal form:

$$
u_1=\frac{u}{w}=x_1^ay^bz_1^{a-1},
v_1=\frac{v}{w}=x_1^cy^dz_1^{c-1},
w_1=w=z_1
$$
at the 3-point $p_1$, which maps to a 3-point.

The analysis of $\Phi_{i+1}$ above points (\ref{eqA7}) and above points satisfying (\ref{eqA4}) where $y=z=0$ are local equations of $E$ 
(and $b,d\ge 1$) is similar. This last case will lead to a form (\ref{eqA5}).
Since $Z$ has  only finitely many components, we inductively construct (\ref{eqA21}).

There now exists a sequence of blow ups of 2-curves $X_r\rightarrow X_m$ which are supported above $q$ such that 
the rational map $X_r\rightarrow Y_1$ is toroidal where ever it is defined, and if ${\cal I}_q{\cal O}_{X_r,p}$ is not invertible, then there
there exist permissible  parameters $x,y,z$ at $p$ for $u,v,w$ such that one of the following forms hold:

$p$ a 1-point 
\begin{equation}\label{eqA14}
u=x^a,
v=x^b(\alpha+y),
w=x^dz
\end{equation}
with $0\ne\alpha\in k$ and $d<\text{min}\{a,b\}$ or

$p$ a 2-point 
\begin{equation}\label{eqA15}
u=x^ay^b,
v=x^cy^d,
w=x^ey^fz
\end{equation}
with $(e,f)<(a,b)<(c,d)$ or $(e,f)<(c,d)<(a,b)$.

We accomplish this as follows. We first consider $u$ and $v$. Suppose that $p\in X_m$ is a 2-point such that ${\cal I}_q{\cal O}_{X_m,p}$ is not invertible. We have forms 
\begin{equation}\label{eqA25}
u=x^ay^b, v=x^cy^d, w=x^ey^fz
\end{equation}
with $e+f>0$ at 2-points $p_i$ above $p$ in the construction of the sequence $X_r\rightarrow X_m$.
At $p_i$ we have an invariant $(a-c)(b-d)$. This is a nonnegative integer if and only if $(a,b)\le(c,d)$ or $(c,d)\le(a,b)$. Further, if $(a-c)(b-d)<0$,
then after blowing up the 2-curve $E$ which has local equations $x=y=0$ at $p_i$, we obtain that all 2-points above $p_i$ have a form (\ref{eqA25}), but $(a-c)(b-d)$ has 
increased. Further $E$ contracts to $q$ on $Y$ since $e+f>0$.

After a finite number of blow ups of 2-curves above $X_m$ (which must contract to $q$) we achieve that all 2-points $p_i$ above a 2-point $p\in X_m$ such that ${\cal I}_q{\cal O}_{X_m,p}$ is not invertible have a form (\ref{eqA25})
with $(a,b)\le(c,d)$ or $(c,d)\le(a,b)$.

We now apply this algorithm to the pairs $u,x^ey^f$ and $v,x^ey^f$ in (\ref{eqA25}) to construct $X_r\rightarrow X_m$.

 We will now inductively construct $X_n\rightarrow X_r$ so that 
${\cal I}_q{\cal O}_{X_n}$ is invertible everywhere and the morphism $X_n\rightarrow Y_1$ is toroidal. We will construct a sequence of blow ups 
\begin{equation}\label{eqA16}
X_n\rightarrow X_{n-1}\rightarrow\cdots\rightarrow X_{r+1}\rightarrow X_r
\end{equation}
so that each $\Phi_{i+1}:X_{i+1}\rightarrow X_i$ is the blow up of a nonsingular curve $\lambda_i$ which is a possible center and is contained in
the locus where ${\cal I}_q{\cal O}_{X_i}$ is not invertible.
We will have that the rational map $f_i:X_i\rightarrow Y_1$ is toroidal where ever it is defined, and all points $p\in X_i$ where ${\cal I}_q{\cal O}_{X_i,p}$
is not invertible have a form (\ref{eqA14}) or (\ref{eqA15}).

Suppose that we have inductively constructed (\ref{eqA16}) up to $X_i$ and ${\cal I}_q{\cal O}_{X_i}$ is not invertible. We will construct $\Phi_{i+1}:X_{i+1}\rightarrow X_i$.

Inspection of the forms (\ref{eqA14}) and (\ref{eqA15}) shows that the locus in $X_i$ where ${\cal I}_q{\cal O}_{X_i}$ is not invertible is a union
of nonsingular curves which are possible centers.  For such a curve $\lambda$, let $\eta$ be a general point of $\lambda$, so that a
form (\ref{eqA14}) holds at $\eta$. Let $A(\lambda)=\text{min}\{a,b\}-d>0$.

Choose a curve $\lambda_i$ which maximizes $A(\lambda)$ on $X_i$. Let $\Phi_{i+1}:X_{i+1}\rightarrow X_i$ be the blow up of $\lambda_i$.
Suppose that $p_i\in\lambda_i$ and $p_{i+1}\in \Phi_{i+1}^{-1}(\lambda_i)$.

First suppose that $p_i$ has the form (\ref{eqA14}). Without loss of generality, we may assume that $a\le b$.  There are   regular parameters $x_1,y,z_1$
in $\hat{\cal O}_{X_{i+1},p_{i+1}}$ satisfying
\begin{equation}\label{eqA17}
x=x_1, z=x_1(z_1+\beta)
\end{equation}
or 
\begin{equation}\label{eqA18}
x=x_1z_1, z=z_1.
\end{equation}

Suppose that (\ref{eqA17}) holds. $p_{i+1}$ is then a 1-point, and 
\begin{equation}\label{eqA19}
u=x_1^a,
v=x_1^b(\alpha+y),
w=x_1^{d+1}(z_1+\beta).
\end{equation}

If $d+1=a=b$ in (\ref{eqA19}), then $X_{i+1}\rightarrow Y_1$ is a morphism near $p_{i+1}$, which maps $p_{i+1}$ to a 1-point, and has a toroidal form
$$
u_1=u=x_1^a,
v_1=\frac{v}{u}-\alpha=y_1,
w_1=\frac{w}{u}-\beta=z_1.
$$

If $d+1=a<b$ in (\ref{eqA19}), then $X_{i+1}\rightarrow Y_1$ is a morphism near $p_{i+1}$, which maps $p_{i+1}$ to a 2-point, and has a toroidal form
$$
u_1=u=x_1^a,
v_1=\frac{v}{u}=x_1^{b-a}(\alpha+y_1),
w_1=\frac{w}{u}-\beta=z_1.
$$

If $d+1<a\le b$ and $\beta\ne 0$ in (\ref{eqA19}) then $X_{i+1}\rightarrow Y_1$ is a morphism near $p_{i+1}$, which maps $p_{i+1}$ to a 3-point, and has a toroidal form obtained from a change of variable in 
$$
u_1=\frac{u}{w}=x_1^{a-d-1}(z_1+\beta)^{-1},
v_1=\frac{v}{w}=x_1^{b-d-1}(\alpha+y_1)(z_1+\beta)^{-1},
w_1=w=x_1^{d+1}(z_1+\beta).
$$

If $d+1<a\le b$ and $\beta= 0$ then (\ref{eqA19}) has a form (\ref{eqA14}) with $d<d+1<\text{min}\{a,b\}$.

Suppose that (\ref{eqA18}) holds. $p_{i+1}$ is then a 2-point, and 
$$
u=x_1^az_1^a,
v=x_1^bz_1^b(\alpha+y),
w=x_1^{d}z_1^{d+1}. 
$$

$X_{i+1}\rightarrow Y_1$ is thus a morphism near $p_{i+1}$, which maps $p_{i+1}$ to a 3-point, and has a toroidal form
$$
u_1=\frac{u}{w}=x_1^{a-d}z_1^{a-d-1},
v_1=\frac{v}{w}=x_1^{b-d}z_1^{b-d-1}(\alpha+y_1),
w_1=w=x_1^{d}z_1^{d+1}.
$$

Now suppose  that $p_i$ has the form (\ref{eqA15}). After possibly interchanging $u$ and $v$, we may assume that $(a,b)< (c,d)$.
After possibly  interchanging $x$ and $y$, we may assume that there are regular parameters $x_1,y,z_1$ in $\hat{\cal O}_{X_{i+1},p_{i+1}}$ satisfying
(\ref{eqA17}) or (\ref{eqA18}) (so that $e<a$).

Suppose that (\ref{eqA17}) holds. Then  $p_{i+1}$ is  a 2-point. We have   
\begin{equation}\label{eqA20}
u=x_1^ay^b,
v=x_1^cy^d,
w=x_1^{e+1}y^f(z_1+\beta).
\end{equation}

If $(e+1,f)=(a,b)$ in (\ref{eqA20}), then $X_{i+1}\rightarrow Y_1$ is a morphism near $p_{i+1}$, which maps $p_{i+1}$ to a 2-point, and has a toroidal form
$$
u_1=u=x_1^ay^b,
v_1=\frac{v}{u}=x_1^{c-a}y_1^{d-b},
w_1=\frac{w}{u}-\beta=z_1.
$$

If $(e+1,f)<(a,b)$ and $\beta\ne 0$ in (\ref{eqA20}), then $X_{i+1}\rightarrow Y_1$ is a morphism near $p_{i+1}$, which maps $p_{i+1}$ to a 3-point, and has a toroidal form obtained from a change of variable in 
$$
u_1=\frac{u}{w}=x_1^{a-e-1}y^{b-f}(z_1+\beta)^{-1},
v_1=\frac{v}{w}=x_1^{c-e-1}y_1^{d-f}(z_1+\beta)^{-1},
w_1=w=x_1^{e+1}y^f(z_1+\beta).
$$

If $(e+1,f)<(a,b)$ and $\beta= 0$ then (\ref{eqA20}) has a form (\ref{eqA15}) with $(e,f)<(e+1,f)<(a,b)$.

Suppose that (\ref{eqA18}) holds. $p_{i+1}$ is then a 3-point, and 
$$
u=x_1^ay^bz_1^a,
v=x_1^cy^dz_1^c,
w=x_1^{e}y^fz_1^{e+1}.
$$

$X_{i+1}\rightarrow Y_1$ is thus a morphism near $p_{i+1}$, which maps $p_{i+1}$ to a 3-point, and has a toroidal form
$$
u_1=\frac{u}{w}=x_1^{a-e}y^{b-f}z_1^{a-e-1},
v_1=\frac{v}{w}=x_1^{c-e}y^{d-f}z_1^{c-e-1},
w_1=w=x_1^{e}y^fz_1^{e+1}.
$$

In summary, we have that all points where ${\cal I}_q{\cal O}_{X_{i+1}}$ is not invertible have a form (\ref{eqA14}) or (\ref{eqA15}) and
if $\lambda_{i+1}\subset \Phi_i^{-1}(\lambda_i)$ is a curve such that ${\cal I}_q{\cal O}_{X_{i+1}}$ is not invertible along $\lambda_i$, we have $0<A(\lambda_{i+1})<A(\lambda_i)$.
Thus after a finite number of blow ups, we construct the desired sequence (\ref{eqA16}), completing the proof of the lemma.
\end{pf}

\begin{Lemma}\label{LemmaB} Suppose that $f:X\rightarrow Y$ is a dominant morphism of nonsingular 3-folds
with toroidal structures determined by SNC divisors $D_Y$ and $D_X=f^{-1}(D_Y)$ such that  $D_X$ contains the singular locus of $f$.
Further suppose
 that $f:X\rightarrow Y$ is toroidal and $C\subset Y$ is a possible center for $D_Y$ which contains a 1-point. Let $\Psi:Y_1\rightarrow Y$ be the blow
up of $C$. Then there exists a commutative
diagram of morphisms
$$
\begin{array}{lll}
\overline X_1&\stackrel{f_1}{\rightarrow}&Y_1\\
\Phi\downarrow&&\downarrow\Psi\\
X&\stackrel{f}{\rightarrow}&Y
\end{array}
$$
such that $\Phi$ is a sequence of possible blow ups for the preimage of $D_X$ supported above $C$ and $f_1$ is toroidal with respect to $D_{Y_1}=\Psi^{-1}(D_Y)$
and $D_{\overline X_1}=\Phi^{-1}(D_X)$.

Further, $\Phi$ has a factorization
$$
\overline X_1=X_n\rightarrow  X_{n-1}\rightarrow \cdots\rightarrow  X_1\rightarrow X
$$
where each $\Phi_{i+1}:X_{i+1}\rightarrow  X_i$ is either the blow up of a section $E_i$ over $C$ such that  ${\cal I}_C{\cal O}_{X_i}$
is not invertible, or $\Phi_{i+1}:X_{i+1}\rightarrow  X_i$ is the blow up of a curve $E_i$ which maps to a 2-point of $Y$ and such that
$E_i$ is contained in the locus where ${\cal I}_C{\cal O}_{X_i}$ is invertible.
\end{Lemma}

\begin{pf} We follow the algorithm of Lemma 18.17 \cite{C3} to construct $\Phi$.

Suppose that $q\in C$ and $p\in f^{-1}(q)$. Then there are permissible parameters $u,v,w$ for $D_Y$ in ${\cal O}_{Y,q}$ and 
regular parameters $x,y,z$ in $\hat{\cal O}_{X,p}$ such that one of the following cases holds:

$q$ is a 2-point and $p$ is a 2-point, 
\begin{equation}\label{eqB1} 
u=x^ay^b, v=x^dy^e,  w=z
\end{equation}
where $uv=0$ is a local equation of $D_Y$ and $u=w=0$ is a local equation of $C$.

$q$ is a 2-point and $p$ is a 1-point, 
\begin{equation}\label{eqB2}
u=x^a, v=x^b(y+\alpha), w=z
\end{equation}
where $0\ne\alpha\in k$, $uv=0$ is a local equation of $D_Y$ and $u=w=0$ is a local equation of $C$.

$q$ is a 1-point and $p$ is a 1-point, 
\begin{equation}\label{eqB3}
u=x^a, v=y,w=z
\end{equation}
where $u=0$ is a local equation of $D_Y$ and $u=w=0$ is a local equation of $C$.

We will construct a sequence of morphisms 
\begin{equation}\label{eqB7}
\cdots\rightarrow X_n\stackrel{\Phi_n}{\rightarrow} X_{n-1}\stackrel{\Phi_{n-1}}{\rightarrow}\cdots\rightarrow X_1\stackrel{\Phi_1}{\rightarrow} X
\end{equation}

where each $\Phi_{i+1}$ is the blow up of a nonsingular curve $E_i$ contained in the locus where ${\cal I}_C{\cal O}_{X_{i}}$ is not invertible,
and for each $q\in C$ and $p\in (f\circ\Phi_1\circ\cdots\circ\Phi_{i})^{-1}(q)$ such that ${\cal I}_C{\cal O}_{X_{i},p}$ is not invertible,
there are permissible parameters $u,v,w$ for $D_Y$ in ${\cal O}_{Y,q}$ and 
permissible parameters $x,y,z$ in $\hat{\cal O}_{X,p}$ such that one of the forms (\ref{eqB4}) - (\ref{eqB6}) below hold.

$q$ a 2-point, $p$ a 2-point 
\begin{equation}\label{eqB4}
u=x^ay^b, v=x^dy^e,
w=x^gy^hz
\end{equation}
with $ae-ba\ne 0$, $(g,h)<(a,b)$.

$q$ a 2-point, $p$ a 1-point 
\begin{equation}\label{eqB5}
u=x^a, v=x^b(y+\alpha),
w=x^dz
\end{equation}
with $0\ne\alpha\in k$, $d<a$,

$q$ a 1-point, $p$ a 1-point 
\begin{equation}\label{eqB6}
u=x^a, v=y,
w=x^dz
\end{equation}
with $d<a$.
Further in the locus where the rational map $X_i\rightarrow Y_1$ is a morphism, $X_i\rightarrow Y_1$ is toroidal. 

Observe that the forms (\ref{eqB1}), (\ref{eqB2}) and (\ref{eqB3}) are special cases of (\ref{eqB4}), (\ref{eqB5}) and (\ref{eqB6}) respectively.

The locus of points where ${\cal I}_C{\cal O}_{X_i}$ is not invertible is a union of nonsingular curves which intersect transversally. If $E$ is a
curve in this locus, and $p'\in E$ is a general point, then $u,v,w$ have a form (\ref{eqB5}) or (\ref{eqB6}) at $p'$. In either case, we define an invariant
$$
\Omega(E)=a-d>0.
$$
Let $\Phi_{i+1}:X_{i+1}\rightarrow X_{i}$ be the blow up of a curve $E_i$ such that $\Omega(E_i)$ is maximal. Suppose that $p_1\in E_i$, $p_2\in \Phi_{i+1}^{-1}(p_1)$
and $q=(f\circ\Phi_1\circ\cdots\circ\Phi_{i})(p_1)$.

Suppose that $p_1$ has a form (\ref{eqB5}). Then $\hat{\cal O}_{X_{i+1},p_2}$ has regular parameters $x_1,y,z_1$ such that 
\begin{equation}\label{eqB8}
x=x_1, z=x_1(z_1+\beta)
\end{equation}
with $\beta\in k$ or 
\begin{equation}\label{eqB9}
x=x_1z_1, z=z_1.
\end{equation}

Suppose that (\ref{eqB8}) holds. Then $p_2$ is a 1-point, 
\begin{equation}\label{eqB10}
u=x_1^a, v=x_1^b(y+\alpha), w=x_1^{d+1}(z_1+\beta).
\end{equation}
If $d+1=a$ in (\ref{eqB10}), then $X_{i+1}\rightarrow Y_1$ is a morphism near $p_2$, mapping $p_2$ to a 2-point, and at $p_2$,  we have a toroidal form
$$
u_1=u=x_1^a,
v=x_1^b(y+\alpha),
w_1=\frac{w}{u}-\beta=z_1.
$$

If $d+1<a$ and $\beta\ne 0$ in (\ref{eqB10}) then 
$X_{i+1}\rightarrow Y_1$ is a morphism near $p_2$, mapping $p_2$ to a 3-point, and at $p_2$,
 we have a toroidal form obtained from a change of variable in
$$
u_1=\frac{u}{w}=x_1^{a-d-1}(z_1+\beta)^{-1},
v=x_1^b(y_1+\alpha),
w_1=w=x_1^{d+1}(z_1+\beta).
$$

If $d+1<a$ and $\beta=0$ in (\ref{eqB10}), then we have a form (\ref{eqB5}) with $d$ increased to $d+1$. The curve $E'$
containing $p_2$ in the locus where ${\cal I}_C{\cal O}_{X_{i+1}}$ is not invertible satisfies
$$
0<\Omega(E')=a-(d+1)<\Omega(E).
$$

Suppose that (\ref{eqB9}) holds. Then $p_2$ is a 2-point. 
\begin{equation}\label{eqB11}
u=x_1^az_1^a,
v=x_1^bz_1^b(y+\alpha),
w=x_1^dz_1^{d+1}.
\end{equation}

Further, 
$X_{i+1}\rightarrow Y_1$ is a morphism near $p_2$, mapping $p_2$ to a 3-point, and at $p_2$,
 we have a toroidal form 
$$
u_1=\frac{u}{w}=x_1^{a-d}z_1^{a-d-1},
v=x_1^bz_1^b(y+\alpha),
w_1=w=x_1^dz_1^{d+1}.
$$
There is a similar argument if $p_1$ satisfies (\ref{eqB6}).

Suppose that $p_1$ has a form (\ref{eqB4}) and $x=z=0$ are local equations of $E_i$ (so that $g<a$). $\hat{\cal O}_{X_{i+1},p_2}$ has regular parameters
$x_1,y_1,z_1$
satisfying (\ref{eqB8}) or (\ref{eqB9}). 

Suppose that (\ref{eqB8}) holds. Then $p_2$ is a 2-point, 
\begin{equation}\label{eqB12}
u=x_1^ay^b,
v=x_1^dy^e,
w=x_1^{g+1}y^h(z_1+\beta).
\end{equation}

If $(g+1,h)=(a,b)$ in (\ref{eqB12}), then 
$X_{i+1}\rightarrow Y_1$ is a morphism near $p_2$, mapping $p_2$ to a 2-point, and at $p_2$,
 we have a toroidal form 
$$
u_1=u=x_1^ay^b,
v_1=v=x_1^dy^e,
w_1=\frac{w}{u}-\beta=z_1.
$$

If $(g+1,h)<(a,b)$ and $\beta\ne 0$ in (\ref{eqB12}), then $X_{i+1}\rightarrow Y_1$ is a morphism near $p_2$, mapping $p_2$
 to a 3-point, and we have a toroidal form obtained from a change of variable in 
$$
u_1=\frac{u}{w}=x_1^{a-g-1}y^{b-h}(z_1+\beta)^{-1},
v=x_1^dy^e,
w_1=w=x_1^{g+1}y^h(z_1+\beta).
$$

If $(g+1,h)<(a,b)$ and $\beta=0$ in (\ref{eqB12}), then (\ref{eqB12}) has the form (\ref{eqB4}) with $g$ increased to $g+1$.

Suppose that (\ref{eqB9}) holds. Then $p_2$ is a 3-point,
$$
u=x_1^ay^bz_1^a,
v=x_1^dy^ez_1^d,
w=x_1^gy_1^hz_1^{g+1}.
$$
$X_{i+1}\rightarrow Y_1$ is a morphism near $p_2$, mapping $p_2$ to a 3-point, and at $p_2$,
 we have a toroidal form 
$$
u_1=\frac{u}{w}=x_1^{a-g}y_1^{b-h}z_1^{a-g-1},
v=x_1^dy^ez_1^d,
w_1=w=x_1^gy_1^hz_1^{g+1}.
$$

By descending induction on $\text{max}(\Omega(E))$, we see that the sequence (\ref{eqB7}) must terminate after a finite number of blow ups,
and we complete the proof of the lemma.
\end{pf}

\section{Preparation}

In this section we prove Theorem \ref{Theorem2}.

We may  assume (after possibly blowing up points on $X$)
that $D_X$ is strongly cuspidal.

To prove this theorem, we may assume by Lemmas \ref{Lemma1} and \ref{Lemma2} that $f$ is prepared (of type 2 a) of Definition \ref{Def31}) above 2 points and 
toroidal above 3-points of $Y$. 
By Corollary \ref{Cor3},
$f$ only fails to be prepared above a finite set of 1-points $\Sigma\subset Y$. Since this reduction involves only blow ups of 2-curves we continue to have the condition that $D_X$ is strongly cuspidal.

Suppose that $q\in \Sigma$. Let $D$ be the component of $D_Y$ containing $q$.
There exists  a very ample effective divisor $L$ on $Y$ such that $q\not\in L$ and $D+L\sim H$ where $H$ is a very ample effective divisor such that $q\not \in H$.
Let $\alpha:Z\rightarrow Y$ be the blow up of $q$, with exceptional divisor $E$. We may replace $L$ with a high multiple of $L$ so that
$\alpha^*H-E$ is very ample on $Z$. Let $N$ be a general member of $\alpha^*H-E$. By Bertini's theorem, $N$ is nonsingular, makes SNCs with 
$D_Z=\alpha^{-1}(D_Y)$, intersects 2-curves of $D_Z$ transversally at general points, does not contain a component of the strict transform on $Z$
of the fundamental locus of $f$, and is disjoint from $\alpha^{-1}(\Sigma-\{q\})$. Let $M=\alpha(N)$. Then $M\sim H$,
 $M$ is nonsingular and  intersects $D$ transversally in a nonsingular curve $\overline \gamma$
which contains $q$, $M$ contains no other points of $\Sigma$, contains no 3-points of $D_Y$, intersects 2-curves of $D_Y$ transversally at general points, does
not contain 
a component of the fundamental locus of $f$ and by Bertini's theorem, at points which are not above $q$, $f^*(M)$ is nonsingular and $f^*(M)+D_X$
is a SNC dvisor. 
After possibly replacing $L$ and $H$ with effective divisors linearly equivalent to $L$ and $H$ respectively, we may assume that
$\overline\gamma\cap(L+H)$ consists of 1-points of $D_Y$ and is disjoint from the fundamental locus of $f$. 

$U=Y-(L+H)$ is an affine neighborhood of $q$. Let $\gamma=\overline\gamma\cap U$. There exist $\overline f,\overline g\in \Gamma(Y,{\cal O}_Y(H))$ such that
$(\overline f)=D+L-H$ and $(\overline g)=M-H$. We can thus define a morphism $\pi:U\rightarrow S={\bf A}^2$ by $\pi(a)=(\overline f(a),\overline g(a))$ for $a\in U$.
Let $\overline q=\pi(q)$. $\pi^{-1}(\overline q) =\gamma$ (scheme theoretically) so $\pi$ is smooth in a neighborhood of $\gamma$. 
 We may thus replace $U$ with an open neighborhood of $\gamma$
so that $\pi$ is smooth.

Let $\overline X=f^{-1}(U)$, and $\overline f=f\mid \overline X$. 
Let $D_U=D_Y\cap U$, $D_U^*=D\cap U$, $D_S^*=\pi(D_U^*)$, $g=\pi\circ\overline f:\overline X\rightarrow S$,
$D_{\overline X}^*=g^{-1}(D_S^*)$, $D_{\overline X}=D_X\cap \overline X$.   
The map $\pi$ is toroidal with respect to $D_S^*$ and $D_U^*$.

 Let $D_1,\ldots, D_m$ be the   components of $D_Y$ other than $D$ which intersect $\gamma$. Since $\gamma$ intersects these components
transversally, we may assume   then that $\pi\mid D_i\cap U$ is \'etale onto its image for $1\le i\le m$.
We further may assume that $\Sigma\cap U=\{q\}$, and (by Bertini's theorem) for $\overline q'\in D_S^*-\{\overline q\}$, there exist regular parameters $u,w$ at $\overline
q'$ such that $u=0$ is a local equation of $D_S^*$, and if $E$ is the curve $w=0$ on $S$, then $E$ is nonsingular, $D_S^*+E$ is a SNC divisor, $g^{-1}(E)$ is nonsingular, and
$g^{-1}(E)+D_{\overline X}$  is a
SNC divisor on $\overline X$. Thus if $q'\in \pi^{-1}(\overline q')$, there exist permissible parameters $u,v,w$ in ${\cal O}_{U,q'}$ (for $D_U$)
such that if $p\in \overline f^{-1}(q')$ then there exist regular parameters $x,y,z$ in $\hat{\cal O}_{\overline X,p}$  such that 
\begin{equation}\label{eq20}
u=x^ay^b, w=z
\end{equation}
where $u=x^ay^b=0$ is a local equation of $D_{\overline X}$ at $p$ (with $a>0$, $b\ge 0$) if $q'\in D-\cup D_i$ and 
\begin{equation}\label{eq30}
u=x^ay^b, v=x^cy^d\gamma, w=z
\end{equation}
where $\gamma\in\hat{\cal O}_{\overline X,p}$ is a unit and $uv=x^{a+c}y^{b+d}=0$ is a local equation of $D_{\overline X}$ at $p$ if $q'\in D\cap D_i$ for some $i$.

Since $\gamma$ intersects the 2-curves $D_i\cap D\cap U$ of $U$ at general points of the 2-curves, after possibly replacing $U$ with a  smaller open neighborhood of $\gamma$, we have that the intersection of the fundamental locus of $f$ with $U$ is contained in $D\cap U$.

We will now establish that  $g$ is toroidal and prepared with respect to $D_S^*$ and $D_{\overline X}^*$ away from the preimages of finitely many
1-points
$\Omega\subset\gamma$ of $D_U$.

Suppose that $q'\in (D_i-D)\cap U$, $p\in \overline f^{-1}(q')$, and $\overline q'=\pi(q')$, which implies that there exist regular parameters $u,w$ at $\overline q'$,
$u,v,w$ at $q'$ such that $v=0$ is a local equation of $D_i$.
$q'$ is not in the fundamental locus of $\overline f$, and $q'$ is a 1-point of $D_U$, so by Abhyankar's lemma there exist regular parameters $x,y,z$ in $\hat{\cal O}_{\overline X,p}$ such that 
\begin{equation}\label{eq32}
u=x, v=y^b,w=z.
\end{equation}
$g$ is defined by $u=x,w=z$ near $p$, which implies that $g$ is smooth, and thus prepared and toroidal for $D_S^*$ and $D_{\overline X}^*$ at $p$.

Suppose that $q'\in (D-\gamma)\cap U$, $p\in\overline f^{-1}(q')$, $\overline q'=\pi(q')$. Then we have a form (\ref{eq20}) or (\ref{eq30}) at $p$, so that  
 $g$ is prepared and toroidal for $D_S^*$ and $D_{\overline X}^*$ at $p$.

Let $\delta=D\cap D_i\cap U$ for some $1\le i\le m$. Suppose that $q'\in\delta\cap \gamma$ and $p\in \overline f^{-1}(q')$.  Then $\pi(q')=\overline q$. Recall that  $q'$ is a general point of the 2-curve $\delta$. There exist regular parameters $u,w$ in ${\cal O}_{S,\overline q'}$
 such that $u=0$ is a local equation of $D$ on $U$, $w=0$ is a local equation of $M$ on $U$, and there exists $v\in {\cal O}_{U,q'}$ 
  such that $v=0$ is a local equation
of $D_i$ and $u,v,w$ are regular parameters in ${\cal O}_{U,q'}$. By our choice of $M$, $\overline f^*(M)$ is nonsingular and makes SNCs with $D_{\overline X}$ at $p$.
Since $uv=0$ is a local equation of $D_{\overline X}$ at $p$, and $w=0$ is a local equation of $\overline f^*(M)$ at $p$, there exist permissible parameters $x,y,z$ in ${\cal O}_{\overline X,p}$
such that 
$$
u=x^ay^b\gamma_1,
v=x^cy^d\gamma_2,
w=z
$$
with  $\gamma_1,\gamma_2$ units in ${\cal O}_{\overline X,p}$. Thus $g$ is prepared and toroidal for $D_S^*$ and $D_{\overline X}^*$ at $p$.

Suppose that $q'\in \gamma$ is a general point. Then $q'$ is a 1-point of $D_U$ and $\overline f$ is finite above $q'$. There exist regular parameters
$u,v,w$ in ${\cal O}_{U,q'}$ such that $u,w$ are permissible parameters for $D_S^*$ at $\overline q=\pi(q')$, and if $p\in \overline f^{-1}(q')$, then
there exist permissible parameters $x,y,z$ in $\hat{\cal O}_{\overline X,p}$ such that 
$$
u=x^a,
v=y,
w=z
$$
by Abhyankar's lemma, which implies that $g$ is prepared and toroidal at $p$ for $D_S^*$ and $D_{\overline X}^*$. 

We conclude that  $g$ is toroidal and prepared with respect to $D_S^*$ and $D_{\overline X}^*$ away from points above finitely many
1-points
$\Omega\subset\gamma$ of $D_U$.

Recall that there are no 3-points of $\overline X$ supported above $D_i\cap U$ for $1\le i\le m$. After blowing up points supported above $\Omega$, we obtain that the
irreducible components $F$ of $D_{\overline X}$ which do not contain a 3-point are precisely the components which dominate $D_i$ for some $i$
or dominate a 2-curve $D_i\cap D$ for some $i$, and the 2-curves $C$ of $D_{\overline X}$ which do not contain a 3-point are precisely the 2-curves which dominate
a 2-curve $D_i\cap D$.

 Suppose that $\Lambda:Z\rightarrow U$ is a dominant morphism of 3-folds, and $D_Z$ is a SNC divisor on $U$. 
We will say that $D_Z$ is $U$ cuspidal if all irreducible components $F$ of $D_Z$ which do not contain a 3-point dominate $D_i$ for some $i$,
or dominate $D_i\cap D$ for some $i$, and the 2-curves $C$ of $D_Z$ which do not contain a 3-point dominate a 2-curve $D_i\cap D$.

By Lemma \ref{Lemma4}, there exists a morphism $\Phi_1:\overline X_1\rightarrow \overline X$ such that $\Phi_1$ is a sequence of possible blow ups 
for the preimage of $D_{\overline X}^*$ of points and
nonsingular curves supported above $\Omega$ such that if $g_1=g\circ \Phi_1:\overline X_1\rightarrow S$ and $f_1=\overline f\circ\Phi_1:\overline X_1\rightarrow U$,
then $g_1$ is prepared  for $D_S^*$ and $D_{\overline X_1}^*=\Phi_1^{-1}(D_{\overline X}^*)$ in a neighborhood of  
all components of $D_{\overline X_1}^*$ which do not map to a point of $\Omega$. 

By blowing up points on components of $D_{\overline X_1}^*$ which dominate a point of $\Omega$, we may suppose that 
$D_{\overline X_1}=\Phi_1^{-1}(D_{\overline X})$
is $U$ cuspidal. 

By 1 of Theorem 3.1 \cite{C5}, there exists a morphism $\Phi_2:\overline X_2\rightarrow \overline X_1$
which is a sequence of possible blow ups for the preimage of $D_{\overline X_1}^*$ of points and nonsingular curves supported above $\Omega$, such that $g_2=\pi\circ f_1\circ\Phi_2:\overline X_2\rightarrow S$ is prepared for $D_S^*$ and
$D_{\overline X_2}^*=\Phi_2^{-1}(D_{\overline X_1}^*)$.  Let $f_2=f_1\circ\Phi_2:\overline X_2\rightarrow U$. We further have that $D_{\overline X_2}=\Phi_2^{-1}(D_{\overline X_1})$
is $U$ cuspidal.
 
Now by 2 of Theorem 3.1 \cite{C5}, there exists a commutative diagram
$$
\begin{array}{lll}
\overline X_3&\stackrel{g_3}{\rightarrow}& S_1\\
\overline \Phi_3\downarrow&&\downarrow\lambda_1\\
\overline X_2&\stackrel{g_2}{\rightarrow}&S
\end{array}
$$
such that $\lambda_1$ is a sequence of possible blow ups for the preimage of $D_S^*$ of points supported above $\overline q$, $\overline\Phi_3$ is a sequence of possible blow ups for the preimage of $D_{\overline X_2}^*$ of points and
nonsingular curves supported above $\gamma$, and $g_3$ is toroidal with respect to $D_{S_1}^*=\lambda_1^{-1}(D_S^*)$ and 
$D_{\overline X_3}^*=\overline\Phi_3^{-1}(D_{\overline X_2}^*)$. 
We further have that 
$D_{\overline X_3}^*$ is $U$ cuspidal.
Let $ f_3=f_2\circ\overline \Phi_3:\overline X_3\rightarrow U$.

Consider the commutative diagram
$$
\begin{array}{lllll}
\overline X_3&\stackrel{\overline f_3}{\rightarrow}&\overline Y_1&\stackrel{\pi_1}{\rightarrow}&S_1\\
\overline\Phi\downarrow&&\overline\Psi\downarrow&&\downarrow\lambda_1\\
\overline X&\stackrel{\overline f}{\rightarrow}&U&\stackrel{\pi}{\rightarrow}&S
\end{array}
$$
where $\overline \Phi=\Phi_1\circ\Phi_2\circ\overline\Phi_3$, $\overline Y_1=U\times_SS_1$ and $\overline\Psi:\overline Y_1\rightarrow U$,
$\pi_1:\overline Y_1\rightarrow S_1$ are the natural projections, and $\overline f_3=f_3\times g_3$. $D_{\overline Y_1}^*=\overline\Psi^{-1}(D_U^*)$ and $D_{\overline Y_1}=\overline\Psi^{-1}(D_U)$ are SNC divisors.
 $\overline Y_1$ is nonsingular, and is obtained from $U$ by possible blow ups for the preimage of $D_U$ of sections over $\gamma$.
Since $g_3$ is toroidal with respect to  $D_{S_1}^*$ and $D_{\overline X_3}^*$, $\overline f_3$ is prepared with respect to $D_{\overline Y_1}^*$
and $D_{\overline X_3}^*$. Over a general point of $\gamma$, $\overline f_3$ is
toroidal with respect to $D_{\overline Y_1}^*$ and $D_{\overline X_3}^*$. Also, over a general point of $\gamma$, $\overline \Phi$ is a sequence of possible blow ups for the preimages of $D_{\overline X}^*$  of sections over $\gamma$.

Recall that $D_Y=D+D_1+\cdots+ D_m+G$, where $G$ consists of the components of $D_Y$ disjoint from $U$, and that the $D_i\cap U$ are \'etale over their images in $S$.
Let $\overline D_i$ be the strict transform of $D_i$ on $\overline Y_1$ for $1\le i\le m$.
$$
D_{\overline Y_1}=\overline\Psi^{-1}(D_U)=D_{\overline Y_1}^*+\overline D_1+\cdots+\overline D_m.
$$
Let $D_{\overline X_3}=\overline\Phi^{-1}(D_{\overline X})$.

We will now verify that $D_{\overline X_3}$ is a $U$ cuspidal SNC divisor on $\overline X_3$ and that $\overline f_3$ is prepared for $D_{\overline Y_1}$
and $D_{\overline X_3}$. Since $\overline f_3$ is prepared for $D_{\overline Y_1}^*$ and $D_{\overline X_3}^*$, we need only verify 
that $\overline f_3$ is prepared for $D_{\overline Y_1}$ and $D_{\overline X_3}$ at points
 $p'\in\overline X_3$ such that $q'=\overline f\circ\overline\Phi(p')\in D_i$ for some $i$.
 
 First suppose that $q'\in D_i-\gamma$ for some $i$. Then $\overline\Phi$ and $\overline\Psi$ are isomorphisms near $p'$ and $q'$ respectively.
 Suppose that 
 $q'\not\in D$. Then  we have permissible parameters $v,u,w$ for $D_U$ at $q'$ which have an expression (\ref{eq32}) at $p'$.
 Thus $\overline f_3$ has an expression 3 of Definition \ref{Def31} at $p'$. Suppose that $q'\in D\cap D_i-\gamma$. Then $q'$ is a 2-point of $D_U$,
  so that $\overline f$ is prepared above $q'$ for $D_U$ and $D_{\overline X}$. Thus $\overline f_3$ is prepared above $q'$
  (for
 $D_{\overline Y_1}$ and $D_{\overline X_3}$).

Suppose that $q'=\overline f\circ\overline\Phi(p')\in\gamma\cap D_i$ for some $i$. 
Without loss of generality, we may assume that $D_i=D_1$.
Recall that $q'\in\gamma\cap D_1$ is a general point of the 2-curve $D\cap D_1$,  
$\overline f$ is prepared  above $q'$ and $\overline f^*(M)$ is nonsingular and makes SNCs with $D_{\overline X}$ above $q'$.
Since $q'\in D\cap D_1$ is a general point, there are no 3-points in $\overline f^{-1}(q')$.
Let $D_1'$ be the reduced divisor on $\overline X$ whose components dominate $D_1$. The irreducible components of $D_1'$ are disjoint by Remark \ref{Remark1}.

 There exist permissible parameters $u,v,w$ in ${\cal O}_{U,q'}$ for the two point $q'$
of $D_U$
such that $u=0$ is a local equation of $D$, $v=0$ is a local equation of $D_1$, $w=0$ is a local equation of $M$ on $U$, and $u,w$ are regular parameters in ${\cal O}_{S,\overline q}$
such that if $p=\overline\Phi(p')\in \overline f^{-1}(q')$, then there exist regular parameters $x,y,z$ in $\hat{\cal O}_{\overline X,p}$
such that one of the following prepared forms for $\overline f$ hold at $p$.   $u,w$ are toroidal forms for $D_U$ and $D_{\overline X}$ in all cases.

\begin{enumerate}
\item $p$ is a 1-point of $D_{\overline X}$ 
\begin{equation}\label{eq1}
\begin{array}{ll}
u&=x^a\\
v&=x^b\gamma\\
w&=z
\end{array}
\end{equation}
where $\gamma\in\hat{\cal O}_{\overline X,p}$ is a unit and $x=0$ is a local equation of $D_{\overline X}$.
\item $p$ is a 2-point of $D_{\overline X}$ which is not on  $D_1'$ 
\begin{equation}\label{eq2}
\begin{array}{ll}
u&=x^ay^b\\
v&=x^cy^d\gamma\\
w&=z
\end{array}
\end{equation}
with $a,b>0$, $\gamma\in\hat{\cal O}_{\overline X,p}$ is a unit and $xy=0$ is a local equation of $D_{\overline X}$.
\item $p$ is a 2-point which is  on  $D_1'$  
\begin{equation}\label{eq3}
\begin{array}{ll}
u&=x^a\\
v&=x^by^c\\
w&=z
\end{array}
\end{equation}
where $xy=0$ is a local equation of $D_{\overline X}$ and $y=0$ is a local equation of $D_1'$.
\end{enumerate}

 $\overline\Psi$ is the sequence of monodial transforms induced by a sequence of quadratic transforms,
$$
S_1=\overline S_n\rightarrow \cdots\rightarrow \overline S_0=S.
$$
Each map $\overline S_{j+1}\rightarrow \overline S_{j}$ is the blow up of the ideal sheaf $m_{j}$ of a point $\overline q_{j}$ above $\overline q$.
Let 
\begin{equation}\label{eq14}
\overline Y_1=\tilde Y_n\rightarrow \cdots\rightarrow \tilde Y_0=U
\end{equation}
be the induced factorization of $\overline\Psi$, where $\tilde\Psi_{j+1}:\tilde Y_{j+1}=U\times_S\overline S_{j+1}\rightarrow \tilde Y_{j}=U\times_S\overline S_{j}$,
is the blow up of a curve $C_j$. Let $\overline \pi_j:\tilde Y_j\rightarrow \overline S_j$ be the natural projection.

$\overline\Phi$ is a sequence of morphisms 
$$
\overline X_3=\tilde X_n\rightarrow \cdots\rightarrow \tilde X_0=\overline X_2\rightarrow \overline X.
$$ 
where $\tilde \Phi_{j+1}:\tilde X_{j+1}\rightarrow \tilde X_{j}$ is a principalization of $m_{j}{\cal O}_{\tilde X_{j}}$, with natural morphism
$\tilde f_j:\tilde X_j\rightarrow \tilde Y_j$. Let $D_{\tilde X_j}$, $D_{\tilde Y_j}$ be the respective preimages of $D_U$, and let $D_{\tilde X_j}^*$, $D_{\tilde Y_j}^*$
be the respective preimages of $D_U^*$. Let $D_{\overline S_j}^*$ be the preimage of $D_S^*$ in $\overline S_j$.
The principalizations $\tilde\Phi_j$ are explicitly described in the proof of Theorem 3.1 \cite{C5}. We have a factorization 
\begin{equation}\label{eqD15}
\tilde X_{j+1}=\hat X_{n_j,j}\rightarrow \cdots\rightarrow \hat X_{0,j}=\tilde X_{j}.
\end{equation}
where each $\hat\Phi_{i+1,j}:\hat X_{i+1,j}\rightarrow \hat X_{i,j}$ is the blow up of a single curve or point $E_{ij}$ which is a possible center for the preimage 
$D_{\hat X_{ij}}^*$ of $D_U^*$
on $\hat X_{i,j}$.  If $E_{ij}$ is a curve, then $E_{ij}$ is in the locus where $m_j{\cal O}_{\hat X_{i,j}}$ is not locally principal. If $E_{ij}$
is a point, then $E_{ij}$ is in the support of $m_j{\cal O}_{\hat X_{i,j}}$ and $m_j{\cal O}_{\hat X_{i,j},E_{ij}}$ is locally principal.
  Further, as is shown in the proof of Theorem 3.1 \cite{C5},
$D_{\tilde X_j}^*$ is $U$ cuspidal (this is the reason for the point blow ups).
Let $D_{\hat X_{ij}}$ be the preimage of $D_U$ on $\hat X_{ij}$.

Recall that $\overline X_2\rightarrow \overline X$ is an isomorphism above $q'$.

We will prove that $D_{\overline X_3}$ is a $U$ cuspidal SNC divisor on $\overline X_3$ and $\overline f_3:\overline X_3\rightarrow \overline Y_1$ is prepared for $D_{\overline Y_1}$ and $D_{\overline X_3}$ above $q'$  by induction on $j$ in the morphisms $\tilde f_j:\tilde X_j\rightarrow \tilde Y_j$.

Recall that we have a fixed choice of regular parameters $u=u_0,v,w=w_0$ in ${\cal O}_{U,q'}$, which are permissible parameters for $D_Y$ at 
the 2-point $q'$, and one of the forms (\ref{eq1}) - (\ref{eq3}) holds at all points of $\overline X$ above $q'$.

Suppose  by induction that $D_{\tilde X_j}$ is a $U$ cuspidal SNC divisor, $\tilde f_j:\tilde X_j\rightarrow\tilde Y_j$ is prepared for $D_{\tilde Y_j}$ and $D_{\tilde X_j}$, and
if $q_j\in \tilde Y_j$ and $\tilde \Psi_1\circ\cdots \circ\tilde\Psi_j(q_j)=q'$,  then
\begin{enumerate}
\item[1.] $q_j$ is a 2-point or a 3-point of $D_{\tilde Y_j}$ and there exist regular parameters $u_j,w_j$ in ${\cal O}_{\overline S_j,\overline q_j'}$,
where $\overline \pi_j(q_j)=\overline q_j'$, such that
$u_j,v,w_j$ are permissible parameters for $D_{\tilde Y_j}$ in ${\cal O}_{\tilde Y_j,q_j}$.
\item[2.] If $q_j\in C_j$, then $u_j=w_j=0$  are local equations of $C_j$ at $q_j$.
\item[3.] If $p_j\in \tilde f_j^{-1}(q_j)$, then there exist regular parameters $x_j,y_j,z_j$ in $\hat{\cal O}_{\tilde X_j,p_j}$ such that one of the
following forms hold:
\vskip .2truein
\noindent {\bf Case 1.} $q_j$ is a 2-point of $D_{\tilde Y_j}$, and $u_jv=0$ is a local equation of $D_{\tilde Y_j}$
(so that $\overline\pi_j(q_j)=\overline q_j'$ is a 1-point), 
and $p_j$ is a 1-point of $D_{\tilde X_j}$ with 
\begin{equation}\label{eqD4} 
u_j=x_j^a, v=x_j^b\gamma_j, w_j=z_j
\end{equation}
where $x_j=0$ is a local equation of $D_{\tilde X_j}$, $\gamma_j$ is a unit in $\hat{\cal O}_{\tilde X_j,p_j}$ or 
$p_j$ is a 2-point of $D_{\tilde X_j}$ with 
\begin{equation}\label{eqD3} 
u_j=x_j^ay_j^b, v=x_j^cy_j^d\gamma_j, w_j=z_j
\end{equation}
where $x_jy_j=0$ is a local equation of $D_{\tilde X_j}$, $\gamma_j$ is a unit in $\hat{\cal O}_{\tilde X_j,p_j}$.

\vskip .2truein
\noindent {\bf Case 2.} $q_j$ is a 3-point of $D_{\tilde Y_j}$, and $u_jvw_j=0$ is a local equation of $D_{\tilde Y_j}$
(so that $\overline\pi_j(q_j)=\overline q_j'$ is a 2-point), 
and $p_j$ is a 1-point of $D_{\tilde X_j}$ with 
\begin{equation}\label{eqD7} 
u_j=x_j^a, v=x_j^b\gamma_j, w_j=x_j^c(z_j+\beta)
\end{equation}
where $x_j=0$ is a local equation of $D_{\tilde X_j}$, $\gamma_j$ is a unit in $\hat{\cal O}_{\tilde X_j,p_j}$ and $0\ne\beta\in k$ or 
$p_j$ is a 2-point of $D_{\tilde X_j}$ with 
\begin{equation}\label{eqD6} 
u_j=x_j^az_j^b, v=x_j^cz_j^d\gamma_j, w_j=x_j^ez_j^f
\end{equation}
where $af-be\ne 0$, $x_jz_j=0$ is a local equation of $D_{\tilde X_j}$, $\gamma_j$ is a unit in $\hat{\cal O}_{\tilde X_j,p_j}$ or
$p_j$ is a 2-point of $D_{\tilde X_j}$ with 
\begin{equation}\label{eqD10} 
u_j=(x_j^ay_j^b)^k, v=x_j^dy_j^e\gamma_j, w_j=(x_j^ay_j^b)^t(z_j+\beta)
\end{equation}
where $x_jy_j=0$ is a local equation of $D_{\tilde X_j}$, $\gamma_j$ is a unit in $\hat{\cal O}_{\tilde X_j,p_j}$, $\text{gcd}(a,b)=1$  and $0\ne\beta\in k$ or 
$p_j$ is a 3-point of $D_{\tilde X_j}$ with 
\begin{equation}\label{eqD11} 
u_j=x_j^ay_j^bz_j^c, v=x_j^dy_j^ez_j^f\gamma_j, w_j=x_j^gy_j^hz_j^i
\end{equation}
where $x_jy_jz_j=0$ is a local equation of $D_{\tilde X_j}$, $\gamma_j$ is a unit in $\hat{\cal O}_{\tilde X_j,p_j}$
and
$$
\text{rank}\left(\begin{array}{lll}a&b&c\\g&h&i\end{array}\right)=2.
$$
\end{enumerate}

We will prove that the above statements hold for $\tilde f_{j+1}:\tilde X_{j+1}\rightarrow \tilde Y_{j+1}$.

Suppose that $q_j\in C_j$ is a 2-point (and  $\tilde\Psi_1\circ\cdots\circ\tilde\Psi_j(q_j)=q'$), so that Case 1 holds, and
$p_j\in \tilde f_j^{-1}(q_j)$. 

${\cal I}_{C_j}{\cal O}_{\tilde X_j,p_j}$ is not invertible, and $u_j,w_j$
satisfy (185) \cite{C3} at $p_j$ if (\ref{eqD4}) holds, $u_j,w_j$ satisfy (190) \cite{C3} at $p_j$ if (\ref{eqD3}) holds and $a,b>0$
(so that $p_j$ is a 2-point of $D_{\tilde X_j}^*$),
$u_j,w_j$ satisfy (185) \cite{C3} at $p_j$ if (\ref{eqD3}) holds and $b=0$ (so that $p_j$ is a 1-point of $D_{\tilde X_j}^*$).

The algorithm of Lemma 18.17 \cite{C3} (as modified after (23) in the proof of Theorem 3.1 of \cite{C5}  by adding appropriate point blow ups to ensure that  $D_{\tilde X_{j+1}}^*$ is $U$ cuspidal)
is applied to construct $\tilde\Phi_{j+1}:\tilde X_{j+1}\rightarrow \tilde X_j$ and $\tilde f_{j+1}:\tilde X_{j+1}\rightarrow \tilde Y_{j+1}$
above $q_j$.
 Suppose that $q_{j+1}\in \tilde\Psi_{j+1}^{-1}(q_j)$, and 
   $\overline\pi_{j+1}(q_{j+1})=\overline q_{j+1}'\in\overline S_{j+1}$.
Then there exist regular parameters $u_{j+1},w_{j+1}$ in ${\cal O}_{\overline S_{j+1},\overline q_{j+1}'}$ such that $u_{j+1},v,w_{j+1}$ are regular parameters
in ${\cal O}_{\tilde Y_{j+1},q_{j+1}}$ and one of the following forms hold:
\vskip .2truein
$\overline q_{j+1}'$ is a 1-point of $D_{\overline S_{j+1}}$ 
\begin{equation}\label{eqD1}
u_j=u_{j+1}, w_j=u_{j+1}(w_{j+1}+\alpha)
\end{equation}
with $\alpha\in k$, or $\overline q_{j+1}'$ is a 2-point for $D_{\overline S_{j+1}}$ 
\begin{equation}\label{eqD2}
u_j=u_{j+1}w_{j+1}, w_j=w_{j+1}.
\end{equation}
If (\ref{eqD1}) holds at $\overline q_{j+1}'$ and $p_{j+1}\in \tilde f_{j+1}^{-1}(q_{j+1})$, then an analysis of the algorithm of
Lemma 18.17 \cite{C3} and Theorem 3.1 \cite{C5} shows that $u_{j+1},v,w_{j+1}$ satisfy one of the forms (\ref{eqD4}) or (\ref{eqD3}) at $p_{j+1}$.

If (\ref{eqD2}) holds at $\overline q_{j+1}'$, and $p_{j+1}\in\tilde f_{j+1}^{-1}(q_{j+1})$, then $u_{j+1},v,w_{j+1}$ satisfy one of the forms
(\ref{eqD7}) - (\ref{eqD11}) at $p_{j+1}$.

If $q_{j+1}\in C_{j+1}$, then $u_{j+1}=w_{j+1}=0$ are local equations of $C_{j+1}$.

Now suppose that $q_j\in C_j$ is a 3-point (and $\tilde\Psi_1\circ\cdots\circ\tilde\Psi_j(q_j)=q'$), so that Case 2 holds, and
$p_j\in \tilde f_j^{-1}(q_j)$.

 If ${\cal I}_{C_j}{\cal O}_{\tilde X_j,p_j}$ is not invertible, then 
after possibly interchanging $u_j$ and $w_j$, then we have one of the following fomrs.  $u_j,w_j$ satisfy (187) \cite{C3} at $p_j$ if (\ref{eqD6}) holds and $a,b>0$,
$u_j,w_j$ satisfy (191) \cite{C3} at $p_j$ if (\ref{eqD6}) holds and $b=e=0$ (in both cases, $p_j$ is a 2-point of $D_{\tilde X_j}^*$).
$u_j,w_j$  satisfy (187) or (191) \cite{C3} if (\ref{eqD11}) holds (so that $p_j$ is a 2 point of $D_{\tilde X_j}^*$). $u_j,w_j$ satisfy (193), (194) or (195) \cite{C3} at $p_j$ if (\ref{eqD11}) holds

The algorithm of  Lemma 18.18 \cite{C3}  is then  applied to construct $\tilde\Phi_{j+1}:\tilde X_{j+1}\rightarrow \tilde X_j$
and $\tilde f_{j+1}:\tilde X_{j+1}\rightarrow \tilde Y_{j+1}$ above $q_j$. 
 Suppose that $q_{j+1}\in \tilde\Psi_{j+1}^{-1}(q_j)$, and $\overline\pi(q_{j+1})=\overline q_{j+1}'\in\overline S_{j+1}$.
Then there exist regular parameters $u_{j+1},w_{j+1}$ in ${\cal O}_{\overline S_{j+1},\overline q_{j+1}'}$ such that $u_{j+1},v,w_{j+1}$ are regular parameters
in ${\cal O}_{\tilde Y_{j+1},q_{j+1}}$ and one of the following forms hold:
\vskip .2truein
$\overline q_{j+1}'$ is a 1-point of $D_{\overline S_{j+1}}$ 
\begin{equation}\label{eqD12}
u_j=u_{j+1}, w_j=u_{j+1}(w_{j+1}+\alpha)
\end{equation}
with $0\ne\alpha\in k$, or $\overline q_{j+1}'$ is a 2-point for $D_{\overline S_{j+1}}$ 
\begin{equation}\label{eqD13}
u_j=u_{j+1}, w_j=u_{j+1}w_{j+1},
\end{equation}
or $\overline q_{j+1}'$ is a 2-point for $D_{\overline S_{j+1}}$ 
\begin{equation}\label{eqD14}
u_j=u_{j+1}w_{j+1}, w_j=w_{j+1}.
\end{equation}

If (\ref{eqD12}) holds at $\overline q_{j+1}'$ and $p_{j+1}\in \tilde f_{j+1}^{-1}(q_{j+1})$, then an analysis of the algorithm of
 Lemma 18.18 \cite{C3}  shows that $u_{j+1},v,w_{j+1}$ satisfy one of the forms (\ref{eqD4}) or (\ref{eqD3}) at $p_{j+1}$.

If (\ref{eqD13}) or (\ref{eqD14}) holds at $\overline q_{j+1}'$, and $p_{j+1}\in\tilde f_{j+1}^{-1}(q_{j+1})$, then $u_{j+1},v,w_{j+1}$ satisfy one of the forms
(\ref{eqD7}) - (\ref{eqD11}) at $p_{j+1}$.  

We have shown that $D_{\tilde X_{j+1}}=\tilde f_{j+1}^{-1}(D_{\tilde Y_{j+1}})$ is a SNC divisor above $q_j$ and that $\tilde f_{j+1}$ is prepared
for  $D_{\tilde Y_{j+1}}$ and $D_{\tilde X_{j+1}}$.

We will now verify that $D_{\tilde X_{j+1}}$ is $U$ cuspidal. Since we are assuming that $D_{\tilde X_j}$ is $U$ cuspidal,
and we know that $D_{\tilde X_{j+1}}^*$ is $U$ cuspidal, we need only verify that
every 2-curve of $D_{\tilde X_{j+1}}$ contained in 
a component of $D_{\tilde X_{j+1}}$ which dominates 
 $D_1$  contains a 3-point. We verify this by induction
in the sequence (\ref{eqD15}).

Let $D_1^{ij}$ be the union of components of $D_{\hat X_{ij}}$  which dominate $D_1$. The irreducible components of $D_1^{ij}$ are
disjoint (by Remark \ref{Remark1}).  Assume that every 2-curve of $D_{\hat X_{ij}}$ which is contained in $D_1^{ij}$
contains a 3-point. We will show that  $D_1^{i+1,j}$  also has this property. 
All points blown up in the construction of $\tilde\Phi_{j+1}$ are either 3-points of $D_{\hat X_{ij}}$ or are 2-points which are  disjoint from 
$D_1^{ij}$.  We may thus assume that the center $E_{i,j}$ blown up by $\hat\Phi_{i+1,j}$
is a curve. 

If $E_{i,j}$ contains a 1-point of $D_{\hat X_{i,j}}$, then $E_{i,j}$ intersects $D_1^{i,j}$ transversally at 2-points of $D_{\hat X_{ij}}$, and thus all 2-curves
of $D_{\hat X_{i+1,j}}$ contained in $D_1^{i+1,j}$ contain a 3-point.
Suppose that $E_{i,j}$ is a 2-curve of $D_{\hat X_{i,j}}$ and $\Lambda$ is an irreducible component of $D_1^{i,j}$.
then either $E_{i,j}$ is contained in $\Lambda$, so that $\Lambda$ contains a 3-point of $E_{i+1,j}$, as $D_{\hat X_{ij}}$ is by assumption
$U$ cuspidal, or else $E_{i,j}$ intersects $\Lambda$ transversally at 3-points of $D_{\hat X_{ij}}$. In either case, all 2-curves of 
$D_{\hat X_{i+1,j}}$ contained in $D_1^{i+1,j}$ contain a 3-point.

We conclude that $D_{\tilde X_{j+1}}$ is $U$ cuspidal.

We have thus established that $\overline f_3$ is prepared for $D_{\overline Y_1}$ and $D_{\overline X_3}$, and $D_{\overline X_3}$ is $U$ cuspidal.

Recall that $\tilde\Phi_{i+1}:\tilde X_{i+1}\rightarrow \tilde X_{i}$ is a principalization of
$m_{i}{\cal O}_{\tilde X_{i}}$ which in a neighborhood of a general point of $\gamma$ is a sequence of blow ups of sections over $\gamma$ where $m_{i}{\cal O}_{\tilde X_{i}}$ is not invertible.

 Each
$\tilde\Psi_{i+1}:\tilde Y_{i+1}\rightarrow \tilde Y_{i}$ is the blow up of a curve $C_i$ which is a section over $\gamma$
and is a possible center for $D_{\tilde Y_{i}}$.

We will construct $\Psi_1:Y_1\rightarrow Y$ such that $\Psi_1^{-1}(U)\cong \overline Y_1$,
$\Psi_1\mid\Psi_1^{-1}(U)=\overline\Psi$ and $\Psi_1^{-1}(D_Y)$ is a SNC divisor by constructing a sequence of morphisms 
\begin{equation}\label{eq15}
Y_1=\hat Y_n\stackrel{\hat\Psi_n}{\rightarrow}\hat Y_{n-1}\rightarrow\cdots\stackrel{\hat\Psi_1}{\rightarrow} Y
\end{equation}
where each $\hat\Psi_{i+1}$ is a product of blow ups of possible centers for the preimage of $D_Y$, and $\hat\Psi_{i+1}^{-1}(\tilde Y_{i})\cong \tilde Y_{i+1}$, $\hat\Psi_{i+1}\mid \tilde Y_{i+1}=\tilde\Psi_{i+1}$ for all $i$.

We will inductively construct (\ref{eq15}). Suppose that we have constructed $\hat\Psi_i:\hat Y_i\rightarrow \hat Y_{i-1}$.

Let $\gamma_{i}$ be the Zariski closure of $C_{i}$ in $\hat Y_i$. Then $\gamma_{i+1}$ is a section over $\overline\gamma$, and is thus a nonsingular
curve. We construct $\hat\Psi_{i+1}$ by first blowing up points on (the strict transform of) $\gamma_{i}$ above $\overline\gamma-\gamma$ where (the strict transform of) $\gamma_{i}$ does not make SNCs with (the preimage of)
$D_{\hat Y_{i}}$, and then blowing up the strict transform of $\gamma_{i}$.

$\overline X_2\rightarrow \overline X$ is an isomorphism away from the preimage of $\Omega$. Thus the sequence of blow ups $\overline X_2\rightarrow X$ extends trivially to a morphism
$X_2\rightarrow X$, so that $X_2\rightarrow X$ is an isomorphism away from the preimage of $\Omega$.

Now we construct $\Phi_3:X_3\rightarrow X_2$ such that $\Phi_3^{-1}(\overline X_2)=\overline X_3$, $\Phi_3\mid\overline X_3=\overline\Phi_3$ 
and $\Phi_3^{-1}(D_{X_3})$ is a SNC divisor by 
constructing a sequence of morphisms
$$
X_3=\hat X_n\stackrel{\hat\Phi_n}{\rightarrow}\hat X_{n-1}\rightarrow\cdots\stackrel{\hat\Phi_1}{\rightarrow}\hat X_0=X_2
$$
where $\hat\Phi_i^{-1}(\tilde X_{i-1})\cong\tilde X_i$ and $\hat\Phi_i\mid \tilde X_i=\tilde\Phi_i$ for all $i$, and so that there are morphisms $\hat X_i\rightarrow \hat Y_i$
which  are toroidal (with respect to the preimages of $D_Y$ and $D_X$) over points of $\overline\gamma-\gamma$. This follows from  application of Lemmas \ref{LemmaA} and \ref{LemmaB},
and the fact that the case when $\gamma_i$ is a 2-curve (or a 3-point is blown up) extends directly to a toroidal morphism.

The resulting morphism $\Phi_3$ is an isomorphism away from the preimage of $\overline\gamma$.

We have constructed a diagram
$$
\begin{array}{lll}
X_3&\stackrel{f_3}{\rightarrow}&Y_1\\
\Phi\downarrow&&\downarrow \Psi_1\\
X&\stackrel{f}{\rightarrow}&Y
\end{array}
$$
such that $\Phi$ and $\Psi_1$ are isomorphisms away from the preimage of $\overline\gamma$ and $f_3$ is prepared with respect to $D_{Y_1}=\Psi_1^{-1}(D_Y)$
and $D_{X_3}=\Phi^{-1}(D_X)$  away from
the points in $\Sigma-\{q\}$. 
 Further, all components of $D_{X_3}$ which do not contain a 3-point and all 2-curves of $D_{X_3}$ which do not
contain a 2-point must contract to points of $\overline\gamma-\gamma$.  In particular, $D_{X_3}$ is cuspidal for $f_3$.
By induction on $|\Sigma|$, we repeat this construction to prove Theorem \ref{Theorem2}.

\section{toroidalization}
In this section we prove Theorems \ref{Theorem3} and \ref{Theorem1}.

\vskip .2truein
\noindent{\bf Proof of Theorem \ref{Theorem3}}

By  Theorem \ref{Theorem2},  we can construct a commutative diagram
$$
\begin{array}{rll}
X_1&\stackrel{f_1}{\rightarrow}&Y_1\\
\Phi_1\downarrow&&\downarrow\Psi_1\\
X&\stackrel{f}{\rightarrow}&Y
\end{array}
$$
such that $\Phi_1$ and $\Psi_1$ are products of possible blow ups  such that $f_1$ is prepared
for $D_{Y_1}=\Psi_1^{-1}(D_{Y})$ and $D_{X_1}=\Phi_1^{-1}(D_{X})$
and $D_{X_1}$ is cuspidal for $f_1$.

Now, we conclude the proof as in the proof of Theorem 0.1 of \cite{C5}.

By descending induction on $\tau(X_2)$ (Definition 2.9 \cite{C5}) and by Theorems 7.11 and 8.1 \cite{C5}, there
exists a commutative diagram
$$
\begin{array}{rll}
X_2&\stackrel{f_2}{\rightarrow}&Y_2\\
\Phi_2\downarrow&&\downarrow\Psi_2\\
X_1&\stackrel{f_1}{\rightarrow}&Y_1
\end{array}
$$
such that $\Phi_2$ and $\Psi_2$ are products of possible blow ups, $f_2$ is prepared
for $D_{Y_2}=\Psi_2^{-1}(D_{Y_1})$ and $D_{X_2}=\Phi_2^{-1}(D_{X_1})$, $D_{X_2}$ is cuspidal for
$f_2$ and $\tau_{f_2}(X_2)=-\infty$.

By Theorem 8.2 \cite{C5}, $f_2$ is toroidal, and the conclusions of the theorem follow.
\vskip .2truein
\noindent{\bf Proof of Theorem \ref{Theorem1}}

 By resolution of singularities and resolution of indeterminacy \cite{H} (cf. Section 6.8 \cite{C6}),
and by \cite{M}, there exists a commutative diagram
$$
\begin{array}{rll}
X_1&\stackrel{f_1}{\rightarrow}&Y_1\\
\Phi_1\downarrow&&\downarrow\Psi_1\\
X&\stackrel{f}{\rightarrow}&Y
\end{array}
$$
where $\Phi_1$, $\Psi_1$ are products of blow ups of points and nonsingular curves supported above $D_Y$, such that $X_1$ and $Y_1$
are nonsingular and projective, and $D_{X_1}=\Phi_1^{-1}(D_X)$ and $D_{Y_1}=\Psi_1^{-1}(D_Y)$ are SNC divisors. 
The proof of Theorem \ref{Theorem1} now follows from Theorem \ref{Theorem3}.

%\end{document}

\vskip.5truein
\noindent
Department of Mathematics

\noindent University of Missouri

\noindent Columbia, MO  65211


\begin{thebibliography}{1000000}
\bibitem[Ab1]{Ab1} Abhyankar, S., {\it On the valuations centered in a local domain}, Amer. J. Math. 78 (1956), 321 -- 348.
\bibitem[Ab2]{Ab2} Abhyankar, S., {\it Algebraic Geometry for Scientists and Engineers}, Amer. Math. Soc., 1990. 
\bibitem[AK]{AK}  Abramovich D.,  Karu K., {\it Weak semistable reduction in characteristic 0}, Invent. Math. 139 (2000),
241 -- 273.
\bibitem[AkK]{AkK} Akbulut, S. and King, H., {\it Topology of algebraic sets}, MSRI publications 25, Springer-Verlag, Berlin.

\bibitem[AKMW]{AKMW} Abramovich, D., Karu, K., Matsuki, K. and Wlodarczyk, J., {\it Torification and factorization of
birational maps}, JAMS 15 (2002), 531 -- 572.
\bibitem[AMR]{AMR} Abramovich, D., Matsuki, K., Rashid, S., {\it A note on the factorization theorem of toric birational maps after Morelli and its toroidal extension}, Tohoku Math J. 51 (1999), 489 -- 537, {\it Correction:}  Tohoku Math J. 52 (2000), 629 -- 631.
\bibitem[BrM]{BrM} Bierstone, E. and Millman, P., {\it Canonical desingularization in characteristic zero by blowing up the maximal strata of a local invariant}, Inv. Math 128 (1997), 207 -- 302.
\bibitem[BEV]{BEV} Bravo, A., Encinas, S., Villamayor, O.,{\it A simplified proof of desingularization and applications}, to appear in
Revista Matematica Iberamericana. 

\bibitem[Ch]{Ch} Christensen, C., {\it Strong domination/weak factorization of three dimensional regular local rings},
Journal of the Indian Math. Soc., 45 (1981), 21 -- 47.
\bibitem[C1]{C1} Cutkosky, S.D., {\it Local factorization of birational maps}, Advances in Mathematics 132 (1997), 167 -- 315.
\bibitem[C2]{C2} Cutkosky, S.D., {\it Local monomialization and factorization of morphisms}, Ast\'erisque 260, 1999.
\bibitem[C3]{C3} Cutkosky, S.D., {\it Monomialization of Morphisms from 3-folds to surfaces}, Lecture Notes in Mathematics 1786, Springer-Verlag, Berlin, Heidelberg, New York, 2002. 
\bibitem[C4]{C4} Cutkosky, S.D., {\it Local monomialization of trancendental extensions},
to appear in the Journal of the Fourier Institute.
\bibitem[C5]{C5} Cutkosky, S.D., {\it Toroidalization of birational morphism of projective 3-folds}, preprint, AG/0407258
\bibitem[C6]{C6} Cutkosky, S.D. {\it Resolution of Singularities}, American Mathematical Society, 2004.
\bibitem[CP]{CP} Cutkosky, S.D. and Piltant, O., {\it Monomial resolutions of morphisms of algebraic surfaces}, Comm. in Alg. 28 (2000), 5935 -- 5959.
\bibitem[CS]{CS} Cutkosky, S.D. and Srinivasan, H. {\it Factorizations of matrices and birational maps}, preprint.
\bibitem[D1]{D1} Danilov, V., {\it Birational geometry of toric 3-folds}, Math USSR Izv. 21 (83), 269 -- 280.
\bibitem[EH]{EH} Encinas, S., Hauser, H., {\it Strong resolution of singularities in characteristic zero}, Comment Math. Helv. 77 (2002), 821 -- 845.
\bibitem[E]{E} Ewald, E., {\it Blow ups of smooth toric 3-varieties}, Abh. math. Sem. Univ. Hamburg 57 (1987).

\bibitem[H]{H} Hironaka, H., {\it Resolution of singularities of an algebraic variety over a field of characteristic zero},
Annals of Math, 79 (1964), 109 -- 326.
\bibitem[K]{K} Karu, K., {\it Local strong factorization of toric birational maps},
J. Alg. Geom 14 (2005), 165 -- 175.
\bibitem[KKMS]{KKMS} Kempf, G., Knudsen, F., Mumford, D., Saint-Donat, B., {\it toroidal embeddings I}, LNM 339, Springer Verlag (1973).
\bibitem[Mat]{Mat}  Matsuki, K., {\it Log resolution of surfaces}, to appear in Contemporary Mathematics.
\
\bibitem[M]{M} Moishezon, B.G. {\it On n-dimensional compact varieties with $n$-algebraic independent meromorphic
functions}, Amer. Math. Soc Translations 63 (1967), 51-177.
\bibitem[Mo]{Mo} Morelli, R., {\it The birational geometry of toric varieties}, J. Algebraic Geometry 5
(1996), 751 -- 782.
\bibitem[O]{O} Oda, T., {\it Torus embeddings and applications}, TIFR, Bombay, 1978.
\bibitem[S]{S} Sally, J., {\it Regular overrings of regular local rings}, Trans. Amer. Math. Soc. 171 (1972)
291 -- 300.

\bibitem[Sh]{Sh} Shannon, D.L., {\it Monoidal transforms}, Amer. J. Math 45 (1973), 284 -- 320.
\bibitem[W1]{W1} Wlodarcyzk, J.,  {\it Decomposition of birational toric maps in blowups and blowdowns},
Trans. Amer. Math. Soc. 349 (1997), 373-411.

\bibitem[W2]{W2} Wlodarcyzk, J., {\it Toroidal varieties and the weak factorization theorem},
Inventiones Math. 154 (2003), 223 -- 331.
\bibitem[Z]{Z} Zariski, O., {\it The compactness of the Riemann manifold of an abstract field of algebraic functions}, Bull. Amer. Math. Soc., 45
(1044), 683 -- 691.
\bibitem[Z1]{Z1} Zariski, O., {\it Introduction to the problem of minimal models in the theory of algebraic surfaces},
Publications of the Math. Soc. of Japan, 1958.
\bibitem[ZS]{ZS} Zariski, O. and Samuel P., {\it Commutative Algebra Volume II}, Van Nostrand, Princeton, 1960.

\end{thebibliography}
\end{document}